\documentclass[11pt]{amsart}
\usepackage{mathrsfs,amsthm,amssymb,amsmath,url,setspace}
\newcommand{\SECTION}[1]{\addtocontents{toc}{\smallskip}\section{#1}\def\SECname{#1}\markboth{#1}{#1}}
\newcommand{\SUBSECTION}[1]{\subsection{#1}\markboth{\SECname -- #1}{\SECname -- #1}}
\newcommand{\SUBSECTIONSTAR}[1]{\subsubsection*{#1}\markboth{\SECname -- #1}{\SECname -- #1}}

\theoremstyle{plain}
    \newtheorem{thm}{Theorem}

    \newtheorem{cor}[thm]{Corollary}
    \newtheorem{fact}[thm]{Fact}
    
    \newtheorem{prob}{Problem}

    \newtheorem{theorem}[thm]{Theorem}
    \newtheorem{Theorem}[thm]{Theorem}
    \newtheorem{Fact}[thm]{Fact}

\theoremstyle{definition}
    
    \newtheorem{nota}[thm]{Notation}
    \newtheorem{Terminology}[thm]{Terminology}
    \newtheorem*{ex}{Example}

    \newtheorem{definition}[thm]{Definition}
    \newtheorem{Definition}[thm]{Definition}

\theoremstyle{remark}

\DeclareMathOperator{\supp}{supp}
\DeclareMathOperator{\pol}{Pol}\DeclareMathOperator{\iinv}{\textbf{Inv}}
\DeclareMathOperator{\Pol}{Pol}\DeclareMathOperator{\ppol}{\textbf{Pol}}

\newcommand{\sm}{\setminus}

\newcommand{\nin}{\notin}

\newcommand{\To}{\rightarrow}

\newcommand{\inv}{^{-1}}

\newcommand{\cl}[1]{\langle #1 \rangle}
\DeclareMathOperator{\Cl}{Cl} 
\newcommand{\CL}{\Cl}

\renewcommand{\O}{{\mathscr O}}
\newcommand{\On}{{\mathscr O}^{(n)}}
\newcommand{\Ok}{{\mathscr O}^{(k)}}
\newcommand{\Oo}{{\mathscr O}^{(1)}}
\newcommand{\Ot}{{\mathscr O}^{(2)}}

\newcommand{\un}{^{(n)}}
\newcommand{\uk}{^{(k)}}
\newcommand{\um}{^{(m)}}
\newcommand{\uo}{^{(1)}}
\newcommand{\ut}{^{(2)}}

\newcommand{\T}{{\mathscr T}}
\newcommand{\C}{{\mathscr C}}
\newcommand{\R}{{\mathscr R}}
\newcommand{\F}{{\mathscr F}}
\newcommand{\J}{{\mathscr J}}

\newcommand{\A}{{\mathscr A}}
\newcommand{\B}{{\mathscr B}}
\newcommand{\D}{{\mathscr D}}
\newcommand{\E}{{\mathscr E}}
\newcommand{\G}{{\mathscr G}}
\newcommand{\M}{{\mathscr M}}

\newcommand{\K}{{\mathscr K}}
\newcommand{\h}{{\mathscr H}}
\renewcommand{\S}{{\mathscr S}}

\renewcommand{\L}{{\frak L}}

\newcommand{\X}{{\frak X}}

\author
{Martin Goldstern}
\address{Algebra\\TU Wien\\Wiedner Hauptstra\ss e 8-10/104\\A-1040 Wien, Austria}
\email{goldstern@tuwien.ac.at}\urladdr{http://www.tuwien.ac.at/goldstern/}
\thanks{The first author is grateful for support through
grant  P17627-N12 of the Austrian Science Foundation (FWF)}
\newcommand{\loor}{\langle \oo 1 \rangle}

\newcommand{\bbn}{{\mathbb N}}
\newcommand{\bbnn}{{\bbn^\bbn}}

\newcommand{\bbr}{{\mathbb R}}

\newcommand{\on}{{\upharpoonright}}

\DeclareMathOperator{\rank}{rank}

\def\itm#1 {\item[{(#1)}]}

\renewcommand{\O}{{\mathscr O}}

\newcommand{\restrictedto}{\upharpoonright}

\newcommand{\ssim}{{\sim}}

\def\oo#1 {\O^{(#1)}}

\def\martinv{\vrule width 3pt height 9pt depth 2pt}
\def\martinhier#1 {\ \endgraf\medskip \leavevmode\martinv \hskip9pt#1 \hskip 9pt \martinv \endgraf\bgroup\sl}
\def\martinweg#1 {\ \endgraf\medskip \leavevmode\martinv \hskip9pt#1 \hskip 9pt \martinv \endgraf\egroup}

\newcommand{\uell}{^{(\ell)}}
\newcommand{\dotminus}{\mathop{%
    \hbox to 0pt{\ \vrule depth -4pt height 5pt width 1pt\hss}-}}

\renewcommand{\dotminus}{\mathop{%
    \raise 2pt\hbox to 0pt{\kern 2.7pt$\cdot$\hss}-}}

\author
{Michael Pinsker}
\address{Algebra\\TU Wien\\Wiedner Hauptstra\ss e 8-10/104\\A-1040 Wien, Austria}
\email{marula@gmx.at}\urladdr{http://dmg.tuwien.ac.at/pinsker/}
\thanks{The second author is grateful for support through
grant P17812 of the Austrian Science Foundation (FWF)} \thanks{This
article is available from \url{http://www.arxiv.org} and from the
authors' websites}

\title[A survey of clones on infinite sets]
    {
        A survey of clones on infinite sets
    }



\begin{document}

    \begin{abstract}
        We summarize what we know about the clone lattice on an
        infinite set and formulate what we consider the most important open problems.
    \end{abstract}

   \maketitle

\tableofcontents \markboth{The clone lattice on an infinite set}{The
clone lattice on an infinite set}

\newpage
\SECTION{The clone lattice on an infinite set}
\label{sec:1}
    A \emph{clone} on a set $X$ is a set of finitary operations on
    $X$ which contains all projections and which is moreover closed under functional composition.
    More formally, define for
    every natural number $n\geq 1$ the set $\On=X^{X^n}$ of all
    operations of arity~$n$, and set $\O=\bigcup_{n\geq 1}\On$ (viewed
    as a disjoint union); then
    $\O$ is the set of all operations of finite arity on~$X$. A
    \emph{projection} is an operation $f(x_1,\ldots,x_n)\in\O$ which
    satisfies an identity of the form $f(x_1,\ldots,x_n)=x_k$ for
    some $1\leq k\leq n$. We denote the $n$-ary projection which projects onto the $k$-th variable by~$\pi^n_k$.
   Whenever $Y$ is a set, $g_1,\ldots, g_n$ are functions from $Y$ to~$X$, and
   $f:X^n\to X$, then we view the tuple $\vec g=(g_1,\ldots, g_n)$ as a
   function from $Y$ to~$X^n$, and write $f(g_1,\ldots, g_n)$ for
   the composition $f\circ \vec g$:
    $$
        f(g_1,\ldots,g_n):\quad
    y\mapsto f(g_1(y),\ldots,g_n(y)).
    $$

    With these definitions, a clone~$\C$ is a
    subset of~$\O$ which contains all $\pi^n_k$ and for which
    the composition $ f(g_1,\ldots,g_n): X^m\to X$
    is an element of~$\C$ whenever $f\in\C$ is $n$-ary and
    $g_1,\ldots,g_n\in\C$ are $m$-ary.
    Examples of clones are:
        \begin{itemize}
        \item The \emph{full clone}~$\O$.
        \item The set $\J$ of all projections on~$X$.
        \item For a partial order $\leq$ on~$X$, the set of
            finitary operations on~$X$ which are monotone with respect
	    to~$\leq$.
        \item More generally, for a relation $R\subseteq X^I$, where $I$ is a not necessarily finite
            index set, the set of all
            operations that preserve this relation. In fact, every clone
            is of this form.
        For infinite~$X$, $I$ can be chosen to be of cardinality $|X|$
        (\cite{Ros71-AClassOfUnivAlgByInfinitaryRel}).
        \item The set of \emph{idempotent} operations (i.e., operations that satisfy the equation
            $f(x,\ldots,x)=x$) on~$X$.
        \item The set of conservative operations on~$X$ (an operation $f\in\On$ is called \emph{conservative}
            iff it satisfies
            $f(x_1,\ldots,x_n)\in\{x_1,\ldots,x_n\}$ for all $x_1,\ldots,x_n\in X$).
        \item For an algebra $\X=(X,\F)$, the set of all endomorphisms, i.e., all operations $g$
            that commute with all operations $f$ of~$\X$: $g(f,\ldots, f) = f(g, \ldots, g)$.
        \item For a topological space $\X=(X,\T)$, if for every $n\geq 1$ we have $X^n$ equipped
            with the product topology, the set of all
            continuous operations which map some product $X^n$ into~$X$.
        \item For an algebra $\X=(X,\F)$, the set of all polynomial
        functions
            of~$\X$.
        \item For an algebra $\X=(X,\F)$, the set of all term
            operations
            of~$\X$. In fact, every clone is of
            this form.
        \end{itemize}
    Observe that the last example is an equivalent definition of a
    clone. Alternatively,
    clones can be defined as precisely the subalgebras of
    a certain algebra $(\O,\{\pi^2_1,\circ,\zeta,\tau,\Delta\})$ with
    base set~$\O$; the latter has been pointed out in slightly
    different form by  
    \textsc{Mal'cev} in \cite{Mal66-IterativeAlgebrasAndPostManifolds}, see also \cite[1.1.1]{PK79}. Since we will
    not make use of this definition, we do not define the
    operations of this algebra here. 

    Ordering all clones on~$X$ by inclusion, one obtains a
    lattice $\Cl(X)$ with largest element $\O$ and smallest element~$\J$. This lattice is called the  \emph{clone
    lattice}. In
    the clone lattice, the meet of two clones is just their intersection,
    while their join is obtained by generating all terms that can be
    built from operations of the two clones. In fact, arbitrary
    intersections of clones are clones again, so the clone lattice
    is complete. Its compact elements are exactly the finitely
    generated clones, so the clone lattice is algebraic since
    clearly every clone is the supremum of its finitely generated
    subclones.

    The clone lattice has only one element if $|X|=1$, but is
    already countably infinite if $|X|=2$. In the latter case,
    $\Cl(X)$ has been completely described by \textsc{Post}
    \cite{Pos41}. If $X$ is finite and has at least three elements,
    then we already have $|\Cl(X)|=2^{\aleph_0}$ \cite{YM59}; it
    cannot be larger since clones are subsets of~$\O$, and since
    $|\O|=\aleph_0$. On infinite~$X$, the size of~$\O$ is~$2^{|X|}$,
    and the clone lattice has cardinality $2^{2^{|X|}}$. This is
    easy to see: Fix two elements $0,1\in X$ and let $f_A\in\Oo$ be the characteristic function of~$A$,
    for every subset~$A$ of~$X$. Then for every set $\A$ of proper non-empty subsets of
    $X\sm\{0,1\}$, consider the clone $\C_\A$ which is generated by
    $\{f_A:A\in\A\}$; a quick check shows that all those clones are
    distinct.

    The larger the base set~$X$, the more complicated the clone
    lattice becomes: If $Y$ is a proper subset of~$X$, then $\Cl(Y)$
    is isomorphic to an interval of~$\Cl(X)$ (see e.g.\ the textbook \cite[Theorem 3.3.5]{PK79} for how to perform such a construction). Even on
    finite base sets with at least three elements, the clone lattice seems to be too complicated
    to be ever fully described. Besides its sheer size, several results indicate that its structure is complex:
    For example, \textsc{Bulatov} has shown that the clone lattice does not satisfy any
    non-trivial lattice identity if $|X|\geq  3$ \cite{Bul93}; it does
    not satisfy any quasi-identity if $|X|\geq 4$ \cite{Bul94}. Also, if
    $|X|\geq 4$, then every countable product of finite lattices
    is a sublattice of~$\Cl(X)$ \cite{Bul94}. Moreover, despite
    considerable effort, fundamental questions such as finding the
    atoms of~$\Cl(X)$ for $|X|\geq 4$ are still open. So far, most attention has
    been given to clones on finite~$X$, since such clones correspond
    to the term operations of finite algebras, which are of particular interest to
    universal algebraists. The three only books on clones,
    the older \cite{PK79} and \cite{Sze86}, and the very recent
    \cite{Lau06}, mainly deal with clones on finite sets. However,
    as can be seen from the examples above, clones appear naturally
    on infinite sets as well, and the study of clones on infinite
    sets allows for the use of not only algebraic methods, but also
    methods from set theory or combinatorial arguments that often
    bear finitary and infinitary aspects. First results on clones on infinite sets date
    back to the 1950s, when \textsc{Yablonskij} studied
    \emph{countable} clones on countably infinite base sets as a
    first generalization from the finite
    (\cite{Jab58-OnLimitLogics} and \cite{Jab59-SomePropertiesOfEnumerableClosedClasses}).
    Recently a considerable number of new
    results has been published. In this survey, we try to summarize
    the state of the art and formulate what we believe are the most
    rewarding open questions in the field.

\SUBSECTION{Conventions and additional notation}
    All notation introduced so far will be valid throughout the paper. In particular, the base set will be~$X$;
    unless mentioned
    otherwise, we will always assume $X$ to be infinite. The sets
    $\On$ (for $n\geq 1$) and $\O$ will sometimes also be denoted by $\On_X$
    and~$\O_X$, respectively, if we want to emphasize the base set~$X$.
    For $\F\subseteq\O$, we write $\cl{\F}$ for the clone
    generated by~$\F$, i.e., the smallest clone that contains~$\F$.
    For~$n\geq 1$, $\F\cap\On$ is called the \emph{$n$-ary fragment}
    of~$\F$ and denoted by~$\F\un$. If $f\in\O$ and the arity of~$f$
    has not yet been 
    given a name, then we denote this arity by~$n_f$. An operation
    $f\in\On$ is called \emph{essentially unary} iff it depends
    only on one of its variables, i.e., iff there is a unary function $F\in\Oo$ and
    $1\leq k\leq n$ such that $f = F( \pi^n_k)$.\\
    For a cardinal~$\lambda$, $[X]^\lambda$ is the set of subsets
    of~$X$ of cardinality~$\lambda$, and $[X]^{<\lambda}$ is the set
    of subsets of~$X$ of cardinality smaller than~$\lambda$. 
    For example, if $X=\bbn $, then $[X]^{|X|}$ is the set of infinite subsets of~$X$, and $[X]^{<|X|}$
    is the set of finite subsets
    of~$X$. We write $\lambda^+$ for the successor cardinal
    of~$\lambda$.

\SUBSECTIONSTAR{Acknowledgement} The authors would like to thank
G\"{u}nther Eigenthaler for his thorough proofreading, and Lutz
Heindorf for his kind help with the Russian literature and his
numerous suggestions.

\SECTION{Non-structure of the clone lattice}
    The size $2^{2^{|X|}}$ of~$\Cl(X)$, the fact that it contains copies of all
    clone lattices over finite sets as intervals, and several
    ``non-structure results'' that will be discussed later in this survey have
    long indicated that the clone lattice is complicated, and that
    there is no hope of describing it completely. Recently, it has been shown that the
    clone lattice is in some sense the most complicated algebraic lattice
    of its size.
    \begin{thm}[Pinsker
    \cite{Pin06AlgebraicSublattices}]\label{thm:pinsker:sublattices:algebraic}
        Let $\L$ be an algebraic lattice which has not more than
        $2^{|X|}$ compact elements. Then $\Cl(X)$ contains a copy of~$\L$ as
        a complete sublattice.
    \end{thm}

    Since the compact elements of~$\Cl(X)$ are just the finitely
    generated clones, their number is not larger than
    $|\O|=2^{|X|}$, and it is easy to construct $2^{|X|}$
    functions that all generate distinct clones: Characteristic
    functions of subsets of~$X$ as used in the introduction are an example. Since it is known that a complete
    sublattice of an algebraic lattice cannot have more compact
    elements than the original lattice, the lattices of Theorem
    \ref{thm:pinsker:sublattices:algebraic} exhaust all complete
    sublattices of the clone lattice. We remark that the analogous
    theorem does not hold on finite~$X$. That is, the clone lattice over finite $X$ has countably infinitely many
    compact elements, but does not have every
    algebraic lattice whose number of compact elements is countable as a
    complete sublattice. For example, the lattice $M_\omega$ (i.e.,
    the lattice consisting of countably infinitely many incomparable
    elements plus a smallest and a greatest element)
    does not embed into the clone lattice over any finite set
    \cite{Bul01}.

    Note that Theorem \ref{thm:pinsker:sublattices:algebraic} implies that $\Cl(X)$ does not satisfy
    any non-trivial lattice (quasi-)identities. This also follows
    from the fact that the latter is the case already for finite $X$ if
    $|X|\geq 3$ ($|X|\geq 4$ for quasi-identities), by a result of
    \textsc{Bulatov} \cite{Bul93} (\cite{Bul94} for
    quasi-identities), since all clone lattices over a finite set
    embed as intervals into~$\Cl(X)$.

    It is open which lattices embed as intervals of the clone
    lattice.

    \begin{prob}\label{prob:algebraiclattice:IsomorphicToInterval}
        Is every algebraic lattice with at most $2^{|X|}$ compact elements
        an interval of the clone lattice~$\Cl(X)$?
    \end{prob}

    If this problem has a positive answer, it would be interesting
    to know whether operations of arity greater than one are
    needed. At this point, we remark that a submonoid of the full
    transformation monoid $\Oo$ can be thought of as a clone whose
    operations are all essentially unary, i.e., depend on at most
    one variable. The lattice of submonoids of~$\Oo$ is therefore
    isomorphic to an interval of the clone lattice, namely the interval
    $[\J,\cl{\Oo}]$ which starts with the projections~$\J$, and
    which has the clone of all essentially unary operations as its
    largest element. Clones in this interval are called
    \emph{unary}.

    \begin{prob}
        Is every algebraic lattice with at most $2^{|X|}$ compact elements
        an interval of the lattice of monoids, i.e., a subinterval of~$[\J,\cl{\Oo}]$?
    \end{prob}

    If the answer to this problem is negative, then it could still
    be the case that Theorem
    \ref{thm:pinsker:sublattices:algebraic} holds for the monoid
    lattice. The proof of the theorem in
    \cite{Pin06AlgebraicSublattices} used non-unary operations.

    \begin{prob}
        Is every algebraic lattice with at most $2^{|X|}$ compact elements
        a complete sublattice of the lattice of monoids, i.e., a subinterval of~$[\J,\cl{\Oo}]$?
    \end{prob}

\SECTION{Precomplete clones}\label{sec:maximal}
    A \emph{precomplete} or \emph{maximal} clone is a dual atom of~$\Cl(X)$, i.e., a clone
    $\C$ that satisfies $\cl{\C\cup\{f\}}=\O$ for all
    operations~$f\nin\C$. 
    A natural question is whether one can determine all precomplete clones
    of the clone lattice. On finite $X$ this question has an
    additional justification: There, since every precomplete clone is determined by its binary fragment (as we will see later in this
    section), the precomplete clones are finite in
    number. Moreover, every clone is contained in a precomplete
    one, which follows immediately from a standard argument using \textsc{Zorn}'s
    lemma and the fact that $\O$ is generated by a finite number of
    functions; in fact, this can be proven easily without the use
    of \textsc{Zorn}'s lemma, see Section \ref{sec:maximal:dual}.
    Therefore, knowledge of the precomplete clones yields an effective
    completeness criterion for finite algebras. That is, given
    an algebra $(X,\{f_1,\ldots,f_k\})$, we can actually decide
    whether every finitary operation from $\O$ is a term operation of the
    algebra: Just check if the operations $f_1,\ldots,f_k$ are all
    contained in one of the precomplete clones. And indeed, the precomplete clones for finite $X$
    have been described in a deep theorem due to  \textsc{Rosenberg} \cite{Ros70}. We will see here that
    on infinite~$X$, the situation is more complicated, but still
    knowledge of some precomplete clones can yield useful completeness
    criteria. The search for such criteria has been pioneered by
    \textsc{Gavrilov} (\cite{Gav59} and \cite{Gav65}) for countably
    infinite~$X$: For example, he defined two binary operations
    $f_1, f_2$, such that an algebra is complete if and only if it
    produces all unary operations and both $f_1$ and $f_2$ as term
    operations.

     One method of describing precomplete clones which was used already
     by \textsc{Gavrilov}, and later also by \textsc{Rosenberg} who
     started investigating precomplete clones on uncountable
     sets in \cite{Ros74}, is
    the following: For $k\geq 1$ and a set of operations
    $\F\subseteq\Ok$, define the clone of \emph{polymorphisms}
    $\pol(\F)$ to consist of all functions $f\in\O$ which satisfy
    $f(g_1,\ldots,g_{n_f})\in\F$ whenever $g_1,\ldots,g_{n_f}\in\F$.
    Then for a clone~$\C$, the following facts (from
    \cite{BKKR69-GaloisTheoryForPostAlgebrasI}, \cite{BKKR69-GaloisTheoryForPostAlgebrasII} for finite $X$ and
    \cite{Ros71-AClassOfUnivAlgByInfinitaryRel} for infinite $X$) are easy to verify:
    $$
        \pol(\C\uk)\supseteq \pol(\C^{(k+1)})\supseteq\C,\quad
        \pol(\C\uk)\uk=\C\uk\quad \text{for all } k\geq 1.
    $$
    In particular,
    $$
        \bigcap_{k\geq 1}\pol(\C\uk)=\C.
    $$
    Therefore, if $\C$ is a precomplete clone such that $\C\uo\neq\Oo$, then
    $\pol(\C\uo)$ is non-trivial and contains~$\C$; whence
    $\C=\pol(\C\uo)$ by the precompleteness. This means that we can find
    all precomplete clones which do not contain all unary operations as
    clones of polymorphisms of some monoid $\M\subseteq\Oo$. If a
    precomplete clone does contain all unary operations, then it does
    not contain all binary operations, for otherwise the well-known fact (see
    \cite{Sie45-SurLesFonctionsDePlusieursVariables})
    that $\cl{\Ot}=\O$ would imply $\C=\O$. Thus by the same
    argument, $\C=\pol(\C\ut)$ and precomplete clones above~$\Oo$ can be
    described as clones of polymorphisms of their binary fragments.

    We will see examples of precomplete clones defined this way later in
    this section. We remark that this representation of
    precomplete clones works also on finite~$X$; hence in that case, there
    exist only finitely many precomplete clones.
    An example of a \emph{collapsing} clone, i.e., of a clone $\C$ for which the descending chain of
    $\pol(\C\uk)$ collapses as a consequence of~$\C=\pol(\C\uo)$, is the one generated by all permutations
    $\S$ of the base set $X$ (see e.g.\ \cite{MP06MinimalAboveS}). On the other hand, we will see in Section \ref{sec:intervals:monoidalIntervals} that
    intervals of the form $[\C,\pol(\C\uo)]$ can also be large.

\SUBSECTION{The number of precomplete
clones}\label{sec:maximal:number}
    Unfortunately, the task of describing all precomplete clones seems rather hopeless. Whereas on finite~$X$,
    the number of precomplete clones is finite, it has been shown
    by \textsc{Gavrilov} in \cite{Gav65} for countable~$X$ and by \textsc{Rosenberg} \cite{Ros76} for uncountable $X$
    that there exist $2^{2^{|X|}}=|\Cl(X)|$ precomplete clones
    on an infinite base set. A shorter and more explicit
    construction proving the latter fact has been provided by \textsc{Goldstern} and
    \textsc{Shelah} in \cite{737}:

    Let $I$ be an ideal of subsets of~$X$,
    that is, a downset of the power set of~$X$ that is closed under
    finite unions. Then the set $\C_I$ of all operations $f\in\O$ which
    satisfy $f[A^{n_f}]\in I$ for all $A\in
    I$ is easily seen to be a clone. Now one can prove that maximal ideals
    (i.e., ideals which cannot be extended to a larger ideal except
    the whole power set of~$X$, or equivalently, ideals dual to ultrafilters) give rise to precomplete clones, and that
    different maximal ideals yield different precomplete clones.
    This immediately implies that the number of precomplete clones is
    $2^{2^{|X|}}$, since it is well-known
    that there exist that many maximal ideals on~$X$.

    \textsc{Cz\'{e}dli} and \textsc{Heindorf} asked which clones~$\C_I$, for $I$
    an ideal, are precomplete if $I$ is not maximal. Define the \emph{support} $\supp(I)$ of an
    ideal $I$ to be the union over the sets of~$I$.
    They found the following criterion for countable~$X$:

    \begin{thm}[Cz\'{e}dli and Heindorf
    \cite{CH01}]\label{thm:czedliheindrf:maximal:ideals}
        Let $I$ be an ideal on a countably infinite base set~$X$.
        \begin{itemize}
            \item
            If $\emptyset\neq \supp(I)\neq X$, then $\C_I$
            is precomplete iff $I$  contains only finite sets or all subsets of~$\supp(I)$.
            \item
            If $\supp(I)= X$, then $\C_I$ is
            precomplete iff $I$ contains some but not all infinite subsets of~$X$
            and for all $B\nin I$ there exists some $f\in \C_I$ such that
            $f[B^{n_f}]=X$.
        \end{itemize}
    \end{thm}

    A drawback of the test for precompleteness in the case where
    the ideal has full support is that in general, one might have to use
    functions of high arity to see that $\C_I$ is precomplete, despite
    the fact that $\C_I$ is actually determined by its unary
    operations:

    \begin{fact}
        For any ideal $I$ on~$X$, $\C_I=\pol(\C_I\uo)$.
    \end{fact}

    Consequently, a solution to the following problem can at
    least be hoped for.

    \begin{prob}
        Find a test for precompleteness of~$\C_I$ that uses unary operations
        only.
    \end{prob}

    If $X$ is uncountable, then the test from
    Theorem~\ref{thm:czedliheindrf:maximal:ideals} does
    not work in general, as can be seen from the following example:

    \begin{ex}
        If $X=\aleph_\omega$, and $I$ is the ideal of bounded
        subsets of~$X$, then one can prove $\C_I$ to be
        precomplete, but for cardinality reasons no countable unbounded subset of
        $X$ can be mapped onto~$X$. For the same reason, the
        precompleteness test fails for $X=\mathbb{R}$ and the ideal
        $J$
        of bounded subsets of~$\mathbb{R}$, although $\C_J$ can be shown to be
        precomplete.
    \end{ex}

    The following problem deals with a well known ideal on the natural numbers.
    \begin{prob}
        Let $X=\mathbb{N}$ be the natural numbers, and define an ideal $I$ to
        consist of all sets $A\subseteq \mathbb{N}$ which have upper
        density~$0$, i.e., for which
        $$
            \overline{\lim}_{n\To\infty}\frac{|A\cap \{0,\ldots , n\}|}{|\{0,\ldots,n\}|}=0.
        $$

        Is $\C_I$ precomplete?
    \end{prob}

\SUBSECTION{Dual atomicity}\label{sec:maximal:dual}

A lattice~$\mathfrak L$ with greatest element~$1$ is called
\emph{dually atomic} iff every element $x\in \mathfrak
L\setminus\{1\}$ is contained in some coatom of~$\mathfrak L$. It is
easy to see that the clone lattice $\CL(X)$ is dually atomic for
finite~$X$. Indeed, for every clone $\C \subsetneqq \O$ there is a
maximal set $\F \subseteq \oo 2 $ with the property $\langle \C \cup
\F\rangle \not= \O$; now $\Pol(\langle \C \cup \F\rangle\ut)$ must
be precomplete.



The above argument used the fact that $\O$ is generated by the finite set
$\O\ut$, which is finite for finite~$X$.
Using \textsc{Zorn}'s lemma one can also easily show the
following for arbitrary $X$:
\begin{Fact}
  Assume that $\F$ is a finite set of functions, $\C_1\subsetneqq
\C_2$ are clones, and $\C_2 = \langle \C_1 \cup \F\rangle$.  Then the
interval $[\C_1,\C_2]$ is dually atomic.
\end{Fact}

For infinite sets $X$ we still have $\O = \langle \O\ut\rangle$, but
since $\O\ut $ is not finite any more, the above arguments cannot be
used, allowing the possibility that the clone lattice for infinite
sets is not dually atomic. It is still not clear whether the
statement
\begin{quote}
 $\CL(X) $ is not dually atomic
\end{quote}
can be proved outright for any infinite set~$X$.  But the following
theorem shows that the above statement cannot be refuted, and that
indeed it holds in many set-theoretic universes, for many sets $X$
(in particular for countably infinite~$X$, assuming the continuum
hypothesis):

\begin{theorem}[Goldstern and Shelah \cite{808}, \cite{884}]
Let $X$ be of regular cardinality~$\kappa$, and assume that
$2^\kappa=\kappa^+$ (in other words: the generalized continuum
hypothesis holds at $\kappa$).
Then the clone lattice $\CL(X)$ is not dually atomic.
\end{theorem}


\begin{prob}
Is it provable (in ZFC, without assumptions on cardinal arithmetic)
that there is a set $X$ whose clone lattice is not dually atomic?
\end{prob}

We also do not know what happens at singular cardinals.
\begin{prob}
Is the clone lattice $\CL(X)$ dually atomic when $|X|$ is singular?
\\
Is it at least consistent that it is (or: is not) dually atomic?
\end{prob}

\SUBSECTION{Precomplete clones that contain all unary
operations}\label{sec:maximal:aboveUnary}

\hyphenation{weak-ly}
\begin{Terminology}
 The binary fragment $\C\cap \Ot $ of a clone~$\C$ is a subset
 of~$\Ot$ containing the two projections and closed
 under the map $(f,g,h) \mapsto f(g,h)$; conversely, every 
 such set is the binary fragment of the clone it generates. 

 Subsets of~$\Ot $ with the above closure properties have been called 
 ``binary Menger algebras'' or ``Menger algebras of rank~2''.   (This name
 has also been used for abstract algebras with a ternary operation 
 $(x,y,z) \mapsto x(y,z)$ satisfying the natural associativity 
 property.)       

 In analogy with the term \emph{monoid}  for unary fragments of 
 clones we suggest the term \emph{dichoid} for binary fragments of 
 clones.     A \emph{binary clone} is a clone generated 
 by its binary fragment; the map $\C\mapsto \C\cap \Ot$ is a natural 
 isomorphism between the lattice of binary clones and the lattice
 of binary Menger algebras (dichoids). 
\end{Terminology}
 


We have mentioned that the precomplete clones $\C$ can be divided
into two classes: Those which do not contain $\oo 1 $ and are
therefore of the form $\Pol(\M) $ for some monoid $\M \subseteq \oo
1 $ (specifically: $\M = \C\uo$), and those which contain $\oo 1 $
and are therefore of the form~$\Pol(\h)$, for some binary Menger
algebra $\h
\subseteq \oo 2 $ (specifically: $\h = \C\ut$). Here we consider
precomplete clones in the second class. We will define them by
describing their binary fragments.

\begin{Fact}\label{fact:maxunary}
\begin{enumerate}
\item
   \label{o1.p}
   Assume that $p\in \oo 2 $  is 1-1.  Then $\langle \oo 1 \cup
   \{p\}\rangle  = \O$.
\item
  \label{o1.max}
   Every clone in~$[\cl{\oo 1} , \O )$  is contained in a precomplete clone.
\item
  \label{o1.clambda}
   For each cardinal $\lambda $ with $2\le \lambda < |X|$, the set
    $$\K_{<\lambda} =  \langle \oo 1 \rangle
          \cup \{ f: |f[X^{n_f}] | <  \lambda \}$$
   is a clone.   For finite numbers $n$ we will write $\K_{n}$ instead
   of~$\K_{<n+1}$.
\end{enumerate}
\end{Fact}
    Observe that (2) follows from (1) using \textsc{Zorn}'s
    lemma. (1) is due to \textsc{Sierpi{\'n}ski}
    \cite{Sie45-SurLesFonctionsDePlusieursVariables}.

 \begin{theorem}[Burle \cite{Bur67}]
  \label{o1.finite}
   If $X$ is finite,
   then the interval $[\langle\oo 1 \rangle , \O ] $  is a chain of
   length $|X|+1$:
$$ \langle \oo 1 \rangle = \K_1 \subsetneqq \B \subsetneqq \K_2 \subsetneqq \ldots \subsetneqq \K_{|X|}=\O$$
(where $\B$ is \textsc{Burle}'s clone described in Section~\ref{sec:local}).

 \end{theorem}

On finite sets we therefore have a unique precomplete clone
above~$\oo 1 $, namely $\K_{|X|-1}$ (this was already
discovered by \textsc{S\l upecki} in \cite{Slu39}; English
translation: \cite{Slu72}). As we will see below, there are two
precomplete clones above $\oo 1 $ when $X$ is countable, and very
many precomplete clones above $\oo 1 $ for most uncountable sets
$X$.

We will now define the two precomplete clones above $\oo 1 $ for a
countable base set~$X$.  While a natural analogue of the first of
these clones can be defined on any set of regular cardinality, the
second clone has an analogue only on so-called \emph{weakly compact}
cardinals.

\subsubsection{Almost unary functions}
To define the first precomplete clone containing~$\Oo$, we consider base
sets $X$ of regular infinite cardinality. We need the following
definition.

\begin{Definition}
 We say that a function $f\in \oo 2 $  is
\emph{almost unary} iff there is a function $F: X \to [X]^{<|X|}$
such that
one of the following holds:
\begin{itemize}
\item For all $x,y\in X$: $f(x,y)\in F(x)$.
\item For all $x,y\in X$: $f(x,y)\in F(y)$.
\end{itemize}
A function $f\in \oo n $ is \emph{almost unary} iff there is a
function~$F$ as above and an index $k\in \{1,\ldots, n\}$ such that
$f(x_1,\ldots, x_n)\in F(x_k)$ for all $x_1,\ldots, x_n\in X$.
\end{Definition}

We assume that $\kappa:= |X|$ is a regular cardinal.  Replacing~$X$ by
$\kappa$ we arrive at the following equivalent definition: A function
$f\in \O_\kappa\un $ is almost unary iff there are $k\in \{1,\ldots,
n\}$ and $F\in\O_\kappa\uo$ such that $f(x_1,\ldots, x_n)\le F(x_k)$
for all $x_1,\ldots, x_n\in\kappa$.

\begin{ex} Any function $f\in \oo 2 _\kappa $ with $f(x,y)=0$ whenever $x<y$ is
almost unary.
\end{ex}

\begin{ex} The function $\min(x,y)\in\Ot_\kappa$  is almost unary, but neither the
function $\max(x,y)\in\Ot_\kappa$ nor the median function
$m_3\in\O^{(3)}_\kappa$ (see Definition~\ref{def:median}) are almost
unary.
\end{ex}

\begin{Definition} We write $T_1 = T_1(X)$ for the set of binary almost
 unary functions.
\end{Definition}

It turns out that the clone of almost unary functions is generated
by its binary fragment (see \cite{Pin04almostUnary}).

\begin{fact}
 The clone $\langle T_1\rangle$ is exactly the clone of all almost unary
functions.
\end{fact}

The following theorem was proved by \textsc{Gavrilov} for
countable~$X$, and by  \textsc{Davies} and  \textsc{Rosenberg} for regular
uncountable~$X$.

\begin{Theorem}[Gavrilov \cite{Gav65}, Davies and Rosenberg \cite{DR85}]
Let~$X$ have regular cardinality.
\begin{enumerate}
\item  $T_1$ is a precomplete binary Menger algebra 
(i.e., a coatom in the lattice
of binary Menger algebras) containing all unary functions.
\item $\Pol( T_1)$ is a precomplete clone containing~$\oo 1 $.
\end{enumerate}
\end{Theorem}

\subsubsection{Never 1-1 functions}

Our next (and, at least on countable sets, last) precomplete clone
 above
$\oo 1 $  can only be defined on base sets of certain
cardinalities.    Again we write $\kappa $ for the cardinality of~$X$, and
then replace~$X$ by~$\kappa$, so we can use a well-order of~$X$.

\begin{Definition}
Let $\kappa$ be a cardinal.   Write \\
\begin{minipage}[t]{8cm}
\begin{itemize}
\item
    $\nabla:= \nabla_\kappa:= \{(x,y)\in \kappa\times \kappa: x<y\}$, the points above the diagonal,
    \item
    $\Delta:= \Delta_\kappa:=\{(x,y)\in \kappa\times \kappa: x>y\}$.
\end{itemize}
\end{minipage}\qquad\qquad
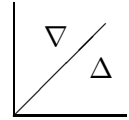
\begin{minipage}[t]{3cm}
\vtop{\ \endgraf
     \setlength{\unitlength}{4144sp}%
     \begingroup\makeatletter\ifx\SetFigFont\undefined%
     \gdef\SetFigFont#1#2#3#4#5{%
       \reset@font\fontsize{#1}{#2pt}%
       \fontfamily{#3}\fontseries{#4}\fontshape{#5}%
       \selectfont}%
     \fi\endgroup%
     \begin{picture}(699,699)(1339,-2848)
     \thinlines
     {\put(1351,-2536){\line( 1, 1){550}}
     }%
     {\put(1351,-1861){\line( 0,-1){675}}
     \put(1351,-2536){\line( 1, 0){675}}
     }%
     \put(1801,-2311){\makebox(0,0)[lb]{$\Delta$}}
     \put(1531,-2121){\makebox(0,0)[lb]{$\nabla$}}
     \end{picture}%
}
\end{minipage}

The notation
 $$\kappa\to (\kappa)^2_n$$ means that for all functions $F:\nabla\to
\{1,\ldots, n\}$ there is a set $A \subseteq \kappa$ of cardinality
$\kappa$ such that $f\on(A\times A)\cap \nabla$ is constant.
\end{Definition}

The fact that $\aleph_0\to (\aleph_0)^2_2$ holds is the statement of
Ramsey's theorem; see \cite{Erdos+Hajnal+Mate+Rado:1984}.
Uncountable cardinals satisfying this partition relation are called
\emph{weakly compact}.

Note that the unary version of this relation is just the pigeonhole
principle: $\kappa\to (\kappa)^1_n$ means that every function from
$\kappa$ to an~$n$-element set is constant on a subset of
$\kappa$ of cardinality~$\kappa$;
 this is true for all infinite cardinals~$\kappa$.

\begin{Definition}
We call a function~$f$ defined on a subset of~$\Delta $ or on a
subset of~$\nabla$ {\em canonical} iff it has one of the following
forms:
\begin{itemize}
\item $f$ is 1-1 (``type 1-1'')
\item $f(x,y) = F(x)$ for some 1-1 function $F\in \O_\kappa\uo$ (``type~$x$'')
\item $f(x,y) = F(y)$ for some 1-1 function~$F\in \O_\kappa\uo$ (``type~$y$'')
\item $f$ is constant (``type $c$'')
\end{itemize}
We call a (partial) binary function $f$ \emph{canonical} if both $f\on \Delta$
and $f\on \nabla$ are canonical.
\end{Definition}
\begin{ex}
The functions $\min$ and $\max$ are  canonical:
  $\min(x,y)=x $ for
$(x,y)\in \nabla$, so $\min\on \nabla$ has type~$x$, whereas $\min\on \Delta $
has type~$y$.
\end{ex}

For the next definition, recall that for a set $Y$ and a cardinal
$\lambda$, $A \in [Y]^{\lambda}$ means that $A$ is a  subset of~$Y$
which has cardinality~$\lambda$.

\begin{fact}\label{fact:canonical}
Let $f\in \O_\kappa\ut$, where $\kappa \to (\kappa)^2_2$.  Then:
\begin{enumerate}
\item There are sets $A,B\in [\kappa]^\kappa$
 such that $f\on (A\times B)$
is canonical.
\item There are 1-1 (strictly increasing) unary functions $u$ and $v$
such that $f(u(x), v(y))$ is canonical.
\end{enumerate}
\end{fact}

We can now classify functions according to the types their canonical
restrictions to various sets of the form $A\times B$ can have. It
turns out that disallowing the type 1-1 will define a binary Menger
algebra.


\begin{Definition}
 We say that a function $f\in \O_\kappa\ut $  is
\emph{never 1-1} iff:
\begin{quote}
Whenever $u,v\in \oo 2 $ are essentially unary functions, \\
then the function~$f(u,v)$, restricted to~$\nabla$, is not 1-1.
\end{quote}
\end{Definition}
Observe that it is sufficient to consider only such pairs $(u,v)$ of
functions where one of the functions depends only on the first
variable, and the other only on the second.

Assuming $\kappa\to (\kappa)^2_2$, we can use
 Fact~\ref{fact:canonical}
to get the following equivalent formulation: $f$ is never 1-1 iff:
for all~$A, B\in [\kappa]^\kappa$ such that $f\on (A\times B)$ is
canonical, the types of~$f\on (A\times B)\cap \Delta$ and $f\on
(A\times B)\cap \nabla$ are among ``type~$x$'', ``type~$y$'' and
``type~$c$'' (but never ``type 1-1'').

More generally, we say that a function $f \in \oo n $ is never 1-1 iff:
\begin{quote}
Whenever $u_1,\ldots, u_n\in \oo 2 $
are essentially unary functions, \\
then the (binary) function $f(u_1,\ldots, u_n)$, restricted to
$\nabla$, is not 1-1.
\end{quote}

We write $\hat T_2$ for the set of all never 1-1 functions, and $T_2$
for the set of all binary never 1-1 functions.

\begin{ex}
Let~$X= \bbn$. The function $p_\Delta\in \oo 2 $, defined by
$$ p_\Delta(x,y) :=  \begin{cases} x^2+y &  \mbox{if $x>y$}\\
0 & \mbox{otherwise,}
   \end{cases}
$$
is canonical; $p_\Delta$ has type 1-1 on~$\Delta$, and type $c$ on
$\nabla$.   Note that $p_\Delta\in T_1 \setminus T_2$.
\end{ex}

\begin{Definition}
We say that a function $f\in \oo 2 _\kappa $  is \emph{densely
unary}, iff:
\begin{quote}
Whenever $A, B\in [\kappa]^{\kappa}$, then there are $A'\in [
A]^{\kappa}$,
 $B'\in [B]^{\kappa}$ such that both $f\restrictedto (A'\times B')\cap \Delta$
and  $f\restrictedto (A'\times B')\cap \nabla $ are essentially unary.
\end{quote}
\end{Definition}

It is clear that every densely unary function is never 1-1.
But for some base sets also the converse holds:
\begin{Theorem}[Gavrilov \cite{Gav65}, Goldstern and Shelah \cite{737}]
Assume that~$X$ is either countable or $|X|$ is a weakly compact cardinal.
 Then:
\begin{itemize}
\item
A  function $f\in \oo 2  $ is densely unary iff $f$ is never 1-1.
\item The set $T_2$ of all binary never 1-1 functions is a binary Menger 
algebra. 
\item The set $\hat T_2$ of all never 1-1 functions is a clone, and
$\hat T_2 = \Pol(T_2)$.
\item $T_2$ is precomplete as a binary Menger algebra.
\item $\Pol( T_2)$ is a precomplete clone.
\item $\Pol( T_1)$ and $\Pol (T_2)$ are the only precomplete clones containing~$\oo 1 $.
\item $ \cl{T_1}$ and $\cl{T_2}$ are the only  coatoms in the lattice of binary clones above~$\oo 1 $.
\end{itemize}
\end{Theorem}
For countable sets this was shown directly by \textsc{Gavrilov} in
 \cite{Gav65}, without using Ramsey's theorem.  \textsc{Goldstern} and
 \textsc{Shelah} in \cite{737} showed this theorem for all cardinals
 satisfying the partition relation $\kappa\to (\kappa)^2_2$.

We remark that
weakly compact cardinals are very large; in particular, they are
regular strong limit cardinals (see Section~\ref{sec:singular}).
  Within the usual framework of Set
Theory (the ZFC axioms) it is not provable that weakly compact
cardinals exist; however, among the ``large cardinals'' considered
in Set Theory, weakly compact cardinals are near the bottom of the
scale, and set theorists working in this area consider the
additional axiom
 ``there are weakly compact cardinals'' (and even much stronger axioms)
plausible (or even ``true'').

\subsubsection{Other regular cardinals}

For uncountable regular cardinals which are not weakly compact, the
situation is unclear. $\Pol( T_1)$ is again a precomplete clone
(assuming that the cardinality of~$X$ is regular), but there seems
to be no reasonable analogue of~$\Pol( T_2)$.  \textsc{Davies} and
\textsc{Rosenberg} proved  in \cite{DR85} that assuming CH the
natural analogue of~$T_2$ is not a binary Menger algebra on the base set
$\aleph_1$, i.e., is not closed under composition.

\begin{Definition}
Let $\lambda$ be a cardinal.  The property $Pr(\lambda)$ is the
following statement: There is a symmetric function $c:\lambda \times
\lambda
 \to \lambda$ with the following
anti-Ramsey property:
\begin{quote}
For all sequences $(a_i:i< \lambda )$ of
pairwise disjoint finite subsets of~$\lambda $, and
for all $c_0\in \lambda$\\
there are $i<j<\lambda $ such that
$c \on (a_i\times a_j) $ is constant with value~$c_0$.
\end{quote}
\end{Definition}

This property is useful for the following theorem:

\begin{theorem}[Goldstern and Shelah \cite{737}]
 Assume $Pr(\lambda)$.  Then there are
$ 2^{ 2^ \lambda} $ many precomplete clones above the unary functions
on the set~$\lambda $.
\end{theorem}
By \cite[III.4 and Appendix 1]{Sh:g},
 $Pr(\lambda)$ holds whenever $\lambda$ is
the successor of an  uncountable regular cardinal (and also for
successors of certain singular cardinals, e.g., if $\lambda=\aleph_{\omega+1}$).

Assuming that weakly compact cardinals exist, the function that
assigns to each cardinal number~$\kappa$ the number of precomplete
clones above ~$\O\uo_\kappa$ is not monotone, as it will often take
the value $2^{2^\kappa}$, but sometimes also the value~$2$.

\subsubsection{Singular cardinals}\label{sec:singular}

\textsc{Davies} and \textsc{Rosenberg} defined in \cite{DR85}
analogues of the clones $T_2$ and $\Pol(T_2)$ on singular strong
limit cardinals.

A cardinal $\kappa$ is called  a \emph{strong limit cardinal} iff
for every cardinal $\lambda <\kappa$ also the cardinality of its
power set $2^\lambda<\kappa$.

\begin{definition}
We say that a function $f $ defined on a product $A\times B $ is
\emph{strictly canonical} iff $f$ has one of the following types:
\begin{itemize}
\item $f$ is 1-1 (``type 1-1'')
\item $f(x,y) = F(x)$ for some 1-1 function $F$ (``type~$x$'')
\item $f(x,y) = F(y)$ for some 1-1 function~$F$ (``type~$y$'')
\item $f$ is constant (``type $c$'')
\end{itemize}
\end{definition}

\begin{nota} For two cardinals $\mu, \nu$ we write $\mu\ll
\nu$ iff $2^{2^{2^\mu}}< \nu$.
\end{nota}

The following fact is a consequence of the
\textsc{Erd\H{o}s}-\textsc{Rado} theorem (see
\cite{Erdos+Hajnal+Mate+Rado:1984}):

\begin{Fact}
Whenever $\mu$ is a cardinal, and $f$ is a function defined on~$A\times
B$, where $\mu^+ \ll |A|=|B|$ then there are sets $A'\in [A]^\mu$,
$B'\in [B]^\mu$ such that $f$ is strictly canonical on~$A'\times B'$.\end{Fact}

As a corollary, we get
\begin{Fact}\label{ER:corollary}
Whenever $|X|$ is a strong limit cardinal and $f\in \O\ut$, then
there are sets $A_\mu, B_\mu\in [X]^\mu$ for unboundedly many
cardinals $\mu<|X|$ such that each restriction $f\on(A_\mu\times
B_\mu)$ is strictly canonical.

Moreover, we may assume that each such canonical restriction has the same type.
\end{Fact}

\begin{definition}
Let $X$ be a set whose cardinality $\kappa$ is a strong limit
cardinal.  We say that $f\in \oo 2 $ is \emph{rarely 1-1} if there is
some $\mu< \kappa$ such that for all $A,B\in [X]^\mu$, the function
$f\on (A\times B)$ is not 1-1.  In other words, we cannot find a
sequence $(A_\mu,B_\mu)_{\mu<\kappa}$ with $|A_\mu| = |B_\mu| = \mu$
such that $f$ is 1-1 on each product $A_\mu\times B_\mu$.

We write $T_2$ for the set of all rarely 1-1 functions.
\end{definition}

\begin{Theorem}[Davies and Rosenberg \cite{DR85}]\label{DR:T2}
  Let $X$ be a set whose cardinality is a strong limit cardinal.
The clone $\cl{T_2}$ is a binary clone containing all unary functions.  If
$|X|$ is moreover singular, then $\Pol(T_2)$ is precomplete.\end{Theorem}
The proof uses Fact~\ref{ER:corollary}.

\begin{prob}
Assume that $X$ has singular cardinality.
How many precomplete clones are there that contain all unary functions?
\end{prob}

\SUBSECTION{Precomplete clones that contain all permutations}
    One step further from determining the precomplete clones above $\Oo$ is
    to describe all precomplete clones that contain the set $\S$ of all
    permutations of~$X$. This amounts to finding all precomplete
    clones above~$\Oo$, which we discussed in the preceding section,
    and those     precomplete
    clones whose unary fragment is a monoid in the interval
    $[\S,\Oo)$ of the monoid lattice. The latter task has been achieved for countably
    infinite $X$ by \textsc{Heindorf} \cite{Hei02} and for
    uncountable $X$ of regular cardinality by \textsc{Pinsker}
    \cite{Pin05MaxAbovePerm}.

    \begin{thm}[Heindorf \cite{Hei02}, Pinsker
    \cite{Pin05MaxAbovePerm}]\label{thm:heindorfpinsker:maximal:aboveS}
        Let $X$ be a set of regular cardinality. The precomplete clones
        which contain all permutations
        but not all unary
        functions are exactly those
        of the form~$\pol(\G)$, where
        $\G\in\{\A,\B,\E,\F\}\cup\{\G_\lambda: 1\leq\lambda\leq|X|,\,\lambda\text{ a cardinal}\}$
        is one of the following submonoids of~$\Oo$:
        \begin{itemize}
            \item{$\A=\{f\in\Oo:|f^{-1}[\{y\}]|<|X|$ for all but fewer than $|X|$ many $y\in X \}$}
            \item{$\B=\{f\in\Oo:|f^{-1}[\{y\}]|<|X|$ for all $y\in X\}$}
            \item{$\E=\{f\in\Oo:|X\sm f[X]|<|X|\}$}
            \item{$\F=\{f\in\Oo: |X\sm f[X]|<|X|$ or $f$ is constant$\}$}
\item{$\G_\lambda=\{f\in\Oo: |X\setminus f[X\setminus A]|\geq \lambda$ for all $A\in
[X]^\lambda \}$}
        \end{itemize}
    \end{thm}

    This shows that there exist relatively few precomplete clones in
    this part of the lattice.

    \begin{cor}\label{cor:pinsker:maximal:aboveS:number}
        Let $X$ have regular cardinality $\aleph_\alpha$.
        Then the number of  precomplete clones on~$X$
        which contain all permutations but not all unary functions is
 $\max(|\alpha|,\aleph_0)$.
    \end{cor}

    Observe that for countably infinite~$X$, as well as for $X$ of
    weakly compact cardinality, Theorem \ref{thm:heindorfpinsker:maximal:aboveS}
    yields a complete description of the clones above
    $\S$ since in that case, the precomplete clones above $\Oo$ are
    known; confer the preceding section. Note also that the
    number of precomplete clones in~$[\cl{\S},\O]\sm [\cl{\Oo},\O]$ is a monotone
    function of the (regular) cardinality of~$X$, whereas the number
    varies in the interval $[\cl{\Oo},\O]$. A metamathematical
    explanation of this could be that the first class of clones is
    determined by unary operations (they are of the form~$\pol(\G)$, for a monoid $\G$), but the second one by binary
    ones (the precomplete clones are of the form~$\pol(\h)$, for $\h\subseteq\Ot$):
    In the first case we make use of the pigeonhole principle, which holds on all infinite~$X$, but in the second case the
    number of precomplete clones depends on binary partition properties of~$|X|$,
    which vary.

    It should be noted here that despite the fact that $\Cl(X)$
    need not be dually atomic, every clone that contains all
    permutations is contained in a precomplete one. This follows as a standard application
    of \textsc{Zorn}'s lemma, because
    $\O$ is finitely generated over~$\S$, confer \cite{Hei02} and
    \cite{Pin05MaxAbovePerm}.

    We know almost nothing about $\Cl(X)$ if the cardinality of the base set is singular. In
    particular, the following is open.

    \begin{prob}
        Generalize Theorem
        \ref{thm:heindorfpinsker:maximal:aboveS} to singular
        cardinals.
    \end{prob}

\SUBSECTION{Symmetric precomplete clones}
    Clones $\C$ which contain all permutations have the property that
    they are \emph{symmetric}, that is, for any permutation $\gamma\in\S$ the
    clone $\C^\gamma$ of all conjugates of functions from~$\C$,
    i.e., of all operations of the form $\gamma\inv
    f(\gamma(x_1),\ldots,\gamma(x_{n_f}))$, where~$f\in\C$, 
    equals~$\C$. Such clones are interesting as they are independent
    of the order, or indeed any other structure that one might
    associate with the base set~$X$.
    \begin{ex}
        If $X$ is the set of natural numbers $\mathbb{N}$, then the clone
        of componentwise (with respect to the natural order) monotone functions is not symmetric and
        does not contain any non-trivial permutation. If $X$ is
        the set of integers $\mathbb{Z}$, then the set of monotone functions is a
        clone which is not symmetric, but which does contain non-trivial permutations.
    \end{ex}
    \begin{ex}
        The clone of idempotent operations is symmetric and does not contain any non-trivial permutation.
    \end{ex}
    \begin{ex}
        If $X$ is countably infinite, then the clone
        $\C_I$ induced by the ideal $I$ as in Section
        \ref{sec:maximal:number} is not symmetric for any ideal $I$ of
        $X$ except when~$\C_I=\O$, which is only the case for 
        the following ideals: The empty ideal, the ideal of finite subsets of~$X$, and the
        ideal of all subsets of~$X$.
    \end{ex}
    \begin{ex}
        The clone of all operations which are either a projection or
        have finite range is symmetric. Also, the clones
        $\K_{<\lambda}$ from Fact \ref{fact:maxunary} are symmetric.
    \end{ex}

    A natural generalization of Theorem
    \ref{thm:heindorfpinsker:maximal:aboveS} is to determine all
    symmetric precomplete clones. \textsc{Pinsker}
    \cite{Pin06PrecompleteConjugation} proved that no new examples occur.

    \begin{thm}[Pinsker \cite{Pin06PrecompleteConjugation}]\label{thm:pinsker:maximal:symmetric}
        Let $X$ be infinite. If $\C\in\Cl(X)$ is symmetric and precomplete, then it contains all
        permutations.
    \end{thm}

    It is readily verified that the symmetric clones form a
    complete sublattice $\Cl_{sym}(X)$ of the clone lattice with
    smallest element $\J$ and largest element~$\O$. We know
    almost nothing about this lattice, except for what follows
    directly from results on~$\Cl(X)$, such as the precomplete clones
    above the permutations. Note that Theorem
    \ref{thm:pinsker:maximal:symmetric} does \emph{not} imply that
    all dual atoms in~$\Cl_{sym}$ contain all permutations, and in
    fact this is not true: The clone of all
    functions $f\in\O$ for which the set $\{x\in X:
    f(x,\ldots,x)\neq x\}$ has fewer than $|X|$ elements is an
    example of a dual atom of~$\Cl_{sym}(X)$ which does not
    contain~$\S$. This has been pointed out in 
    \cite{Pin06PrecompleteConjugation} and is a consequence of the
    complete description of the interval of~$\Cl(X)$ above this
    clone in
    \cite{747}, see Section \ref{sec:aboveIdempotent}.

    \begin{prob}
        Determine the dual atoms of~$\Cl_{sym}(X)$.
    \end{prob}

    We remark that whereas $\Cl(X)$ need not be dually atomic, the
    sublattice $\Cl_{sym}$ is; this is 
    a consequence of \textsc{Zorn}'s lemma and the fact
    that $\O$ is finitely generated in~$\Cl_{sym}(X)$, i.e.,
    there exist finitely many functions such that $\O$ is
    the only symmetric clone that contains those functions \cite{Pin06PrecompleteConjugation}.

    It might be interesting to note that on finite~$X$, all
    symmetric clones are known (\cite{Kho92},\cite{Kho93},\cite{Kho94},\cite{Mar96a},\cite{Mar96b},
    see also the survey paper \cite{Sze03}). If $X$ has at least five elements, then
    the only symmetric precomplete
    clone is the \textsc{S\l upecki}-clone of all
    functions which are either essentially unary or take at most $|X|-1$
    values. In this case the
    clone of all idempotent functions is the only other
    clone which is precomplete in~$\Cl_{sym}(X)$ (but not in~$\Cl(X)$,
    since it is properly contained in
    the clone of all $f\in\O$ for which $f(a,\ldots,a)=a$, for any fixed $a\in X$).
    For $|X|<5$ the situation is more complicated, see
    \cite{Sze03}.

\SUBSECTION{The rank of clones}

    A possible next step after finding a precomplete clone is to try to
    determine its maximal subclones, which one can imagine as the
    clones in the second level from above (with $\O$ being at level
    zero and the precomplete clones being at level one). More generally,
    \textsc{Gavrilov} \cite{Gav74-ThePowerOfTheSetOfClassesOfFiniteHeight} inductively
    defined the \emph{rank} of a clone as follows: $\O$ has rank~$0$,
 and a clone $\C$ is said to be of rank~$n$, for $n\geq 1$
    a natural number, iff for any function $f\nin \C$ the clone
    $\cl{\{f\}\cup\C}$ has rank at most~$n-1$, and there exists at
    least one function $f\nin\C$ such that the rank of
    $\cl{\{f\}\cup\C}$ equals~$n-1$.   Thus, the rank of a clone $\C$
    measures a kind of distance from 
    $\C$ to~$\O$.  Observe that there exist clones which do not have such a
    finite rank: Any countably generated clone is an example.

    \begin{thm}[Gavrilov
    \cite{Gav74-ThePowerOfTheSetOfClassesOfFiniteHeight}]\label{thm:gavrilov:submaximal:finiteHeight}
        Let $X$ be countably infinite. Then for every~$n\geq 1$, there exist $2^{2^{\aleph_0}}$
        clones of rank~$n$.
    \end{thm}

    One can extend this idea and inductively define a clone $\C$ to
    have rank~$\alpha$, where $\alpha$ is an ordinal, iff it does not
    have rank 
    $\beta$ for any ordinal $\beta<\alpha$, but every clone properly extending $\C$
    does have rank $\beta$ for some $\beta <\alpha$. If it does not
    have rank $\alpha$ for any ordinal~$\alpha$, then we define its
    rank to be~$\infty$. It is clear that there exist clones of rank~$\infty$,
    since there exist clones which have infinite
    ascending chains above them. Also, if a clone has rank~$\infty$,
    then all of its subclones have rank~$\infty$. We remark that rank
    $\infty$ is not to be confused with a ``proper'' transfinite
    rank, such as~$\omega$.

    \begin{ex} The clone $\cl{T_1}$ of
    Section~\ref{sec:maximal:aboveUnary} has rank~$\omega$, see Theorem \ref{pinsker:chain}.
    \end{ex}

    It might be interesting to observe that in the clone lattice of the
    two-element set, every clone has a proper rank and $\J$ has
    rank~$\omega+2$; moreover, for each $\alpha\leq \omega+2$ there
    exist only finitely 
    many clones of rank~$\alpha$. If $X$ has at least three elements,
    then we already have clones with rank~$\infty$. 
    Back to infinite~$X$, we remark that
    the rank of a clone (if it is proper)
    must be an ordinal smaller than $(2^{|X|})^+$ (the successor cardinal of~$2^{|X|}$).

    \begin{prob}
        Let $\alpha< (2^{|X|})^+$. How many clones of rank $\alpha$
        are there?  What is the first $\alpha$ such that there are
    no clones of rank~$\alpha$?
    \end{prob}

    If $\D\subseteq\C$ are clones, then analogously to the usual rank of a
    clone, one can define the \emph{rank of~$\D$ in~$\C$}. So given
    a clone $\C$ and an ordinal~$\alpha$, one can ask how many
    subclones of~$\C$ there are which have rank $\alpha$
    in~$\C$. For~$\alpha=1$, this is just the question of how many
    maximal proper 
    subclones $\C$ has.

    \textsc{Marchenkov} studied this relative rank for finite
    ordinals, on a countably infinite base set: Identify $X$ with the
    set of natural numbers $\mathbb{N}$. Call a clone $\C$ on
    $\mathbb{N}$ \emph{elementarily closed} (in the sense of \textsc{Skolem})
    iff it contains the constant function with value~$1$, the addition~$x+y$, the
    operation $x\dotminus y:=\max(x-y,0)$, and together with each
    function $f(x_1,\ldots,x_n,y)\in\C$ also the function
    $g(x_1,\ldots,x_n,y)=\sum_{i=0}^y f(x_1,\ldots,x_n,i)$. Then we
    have:

    \begin{thm}[Marchenkov
    \cite{Mar81-CardinalityOfTheSetOfPrecomplete}]
        Let $X=\mathbb{N}$, and let $\C$ be a clone of cardinality
        $\lambda$ which is elementarily closed. Then for all
        finite~$n\geq 1$, there exist $2^\lambda$ subclones of~$\C$
        which 
        have rank $n$ in~$\C$.
    \end{thm}

    Observe that this theorem implies Theorem
    \ref{thm:gavrilov:submaximal:finiteHeight}; in particular, it
    implies that there exist $2^{2^{\aleph_0}}$ precomplete clones
    on a countably infinite base set.

\SECTION{Minimal clones}

    A clone is called \emph{minimal} iff it is an atom in~$\Cl(X)$,
    i.e., iff the only clone below it is the clone of 
    projections. Clearly every minimal
    clone is generated by a single operation. We call operations
    which generate minimal clones, and which have minimal arity in
    the sense that no operation of smaller arity generates the same
    clone,
    \emph{minimal} as well. An obvious necessary
    and sufficient condition for an operation $f\in\On$ to be minimal is
    that all non-trivial terms which it generates have arity at least $n$ and generate~$f$.

    \begin{ex}
        For every~$n\geq 2$, let $f_n\in\S$ be a permutation which has only one finite cycle of
        length~$n$, and which is the identity otherwise. If $n$ is
        a prime number, then it is easy to see that all non-trivial iterates
        of~$f$ can reproduce~$f$, so $f$ is minimal. If $n$ is divisible by some $k$ with $2\leq
        k<n$, then the function $f^k$ obtained by iterating $f$ $k$
        times cannot generate~$f$, hence $f$ is not minimal.
    \end{ex}

    On finite~$X$, the clone lattice is atomic,
    which is to say that every clone (except for the clone of projections) contains a minimal clone; see
    e.g.\
    the survey papers \cite{Csa05-MinimalClonesAMinicourse} or
    \cite{Qua95-ASurveyOfMinimalClones} for a proof.
    If $X$ is infinite, then this is not so, as can be seen from
    the following simple example:

    \begin{ex}
        Let $f\in\S$ be a permutation which has only infinite
        cycles, so if $f^k$ (where $k\geq 1$) is any iterate of~$f$, then $f^k(x)\neq x$
        for all~$x\in X$. The interval $[\J,\cl{\{f\}}]$ is
        isomorphic to the lattice of all submonoids of the
        monoid $(\mathbb{N},+,0)$. In particular, it is not atomic.
    \end{ex}

    This example also shows in an easy way that $\Cl(X)$ does not
    satisfy any non-trivial lattice identity, since
    the submonoid lattice
    of~$(\mathbb{N},+,0)$ does not (the latter was shown in \cite{RK88-CommutativeSemigroups}).

    Even on finite~$X$, despite the fact that there exist only finitely many minimal clones (see e.g.\ the textbook
    \cite{PK79} or the surveys \cite{Csa05-MinimalClonesAMinicourse}
    and
    \cite{Qua95-ASurveyOfMinimalClones}),
    no explicit list of the minimal
    operations is known. Solving this problem is even more
    difficult on infinite~$X$, since every minimal operation $f$ on a
    finite set $Y$ can be extended to a minimal operation on infinite
    $X$: Just observe that the operation of every algebra in the variety generated by the algebra
    $(Y,\{f\})$ is minimal, since minimality can be read off the
    equations $f$ satisfies (that is, minimality is an abstract property, see Fact~\ref{abstract:minimal}).
    It therefore suffices to take a
    subalgebra of a product of~$(Y,\{f\})$ which has cardinality
    $|X|$ to interpret $f$ on~$X$. Therefore, finding all minimal operations on~$X$ includes
    finding all minimal operations on all finite sets, and
    the following seems very ambitious:

    \begin{prob}
        Describe the minimal clones of~$\Cl(X)$.
    \end{prob}

    Some necessary properties for a function to be
    minimal have been described by \textsc{Rosenberg}
    \cite{Ros83-MinimalClonesI}: For example, it is clear that
    every minimal operation which depends on more than one
    variable must be idempotent, for otherwise it generates a
    non-trivial unary operation, which in turn cannot generate the
    original operation. Call an operation $f\in\On$ a \emph{semiprojection} iff there exists $1\leq k\leq n$ such that
    $f(x_1,\ldots,x_n)=x_k$ whenever $|\{x_1,\ldots,x_n\}|<n$. We say that $f\in\O^{(3)}$ is a \emph{majority operation}
    iff it satisfies the equations
    $f(x,x,y)=f(x,y,x)=f(y,x,x)=x$. The proof in
    \cite{Ros83-MinimalClonesI}, although formulated for finite
    $X$ there,
     works on infinite sets as well, which yields the following theorem:

    \begin{thm}[Rosenberg
    \cite{Ros83-MinimalClonesI}]\label{thm:rosenberg:minimal:fivetypes}
        Every minimal operation is of one of the following types:
        \begin{itemize}
            \item A unary operation $f$ which either satisfies $f^2=f$ or is
            a non-trivial permutation such that $f^p$ is the identity for some prime number~$p$.
            \item A binary idempotent operation.
            \item A ternary majority operation.
            \item A ternary reduct $x+y+z$ of an elementary
            $2$-group~$(X,+)$.
            \item An $n$-ary semiprojection~$(n> 2)$.
        \end{itemize}
    \end{thm}

    \textsc{P\'{a}lfy} \cite{Pal86-TheArityOfMinimalClones} proved the existence of a minimal
    operation of arity $n$ for every $1\leq n\leq|X|$, for finite~$X$;
    since the proof works for infinite $X$ as well, there exist
    minimal operations of all arities here. The latter also follows from
    the fact that on finite~$X$, every non-trivial semiprojection of arity
    $|X|$ generates a minimal operation, which in turn must have arity $|X|$ as well, and from the extension of minimal
    operations to infinite $X$ as described above.

    Of course, on finite $X$ every clone has only finitely many
    binary operations. The following question for infinite $X$ has been posed in slightly different form
    in \cite{Csa05-MinimalClonesAMinicourse}:

    \begin{prob}\label{prob:minimal:infinitelyManyBinary}
        Does there exist a minimal clone with an infinite number of
        binary operations?
    \end{prob}

    \textsc{Machida} and \textsc{Rosenberg} extended the notion of
    a minimal clone and defined an \emph{essentially minimal clone}
    to be a clone which contains at least one operation which is not
    essentially unary, and whose proper subclones do not have this property. Rechecking this
    definition, one sees that a clone is essentially minimal iff it is an
    atom of a monoidal interval (see Section \ref{sec:intervals:monoidalIntervals}
    for information on
    monoidal intervals). In \cite{MR92}, they exhibited an essentially
    minimal clone with an infinite number of binary operations; the unary fragment of this clone is non-trivial.
    In this light, Problem \ref{prob:minimal:infinitelyManyBinary} is the
    question whether the same result can be achieved with a
    trivial monoid (unary fragment).

\SECTION{Intervals} \SUBSECTION{The interval of clones that contain
all unary operations}\label{sec:interval:aboveOo}


In this section we mainly concentrate on a countable base set~$X$,
say $X= \bbn = \{0,1,\ldots\}$ (except for
Theorem~\ref{pinsker:chain}, which
 works on all sets of
 regular cardinality).  Notwithstanding the fact that there are only two
precomplete clones in the interval~$[\loor , \O ]$, there are
several results indicating that this interval is extremely
complicated. Since the clones $\pol(T_1)$ and $\pol(T_2)$ are the
only precomplete clones in the interval, and since the interval is
dually atomic, it can be written as
$$ [\loor , \O] \ \ = \ \ [\loor , \Pol(T_1)] \ \cup \ [\loor ,
\Pol(T_2)]\ \cup \ \{\O\}. $$ Hence an analysis of this interval
naturally splits into two (overlapping) areas: clones below
$\Pol(T_1)$ and clones below~$\Pol(T_2)$.

\subsubsection{Around the binary clone $\cl{T_1}$}

Recall that $T_1$ is the set of all binary functions which are almost
unary. We identify our base set $X$ with the cardinal number~$|X|$, so
that $X$ is linearly ordered.

\begin{Definition}\label{def:median}
For any $\vec x = (x_1,\ldots, x_k)\in X^k$ ($k\ge 2$) we let
 $\sigma_{\vec x}$ be any permutation of~$\{1,\ldots, k\}$ such
 that~$x_{\sigma_{\vec x}(1)} \le \cdots \le x_{\sigma_{\vec x}(k)}$, and we
 write  $m_k$ for the $k$-ary function~$\vec x \mapsto x_{\sigma_{\vec x}(2)}$.
\\
Thus, $m_2(x,y) = \max(x,y)$ and $m_3(x,y,z)$ is the median
of~$x,y,z$. In general, $m_k(x_1,\ldots, x_k)$ is the ``second
smallest'' element of the $k$-element multiset~$\{x_1, \ldots, x_k\}$.\end{Definition}

The following theorem completely describes the clones containing
$T_1$.

\begin{Theorem}[{Pinsker} \cite{Pin04almostUnary}]
\label{pinsker:chain}
Let $|X| = \kappa$ be a regular cardinal.  We will consider clones on
$\kappa$.
\begin{enumerate}
\item
   Let $p_\Delta\in \oo 2 $ be a function which is 1-1 on~$\Delta$,
   and constantly 0 on~$(X \times X) \setminus \Delta$.  Then $T_1$ is
   generated by~$\oo 1 \cup \{p_\Delta\}$.
\item
   The median function $m_3$ is a ternary
   function in~$\Pol(T_1)$ which is not generated by functions
   in~$T_1$. Hence~$\langle T_1\rangle \subsetneqq \Pol(T_1)$.
\item The precomplete clone  $\Pol(T_1)$ is generated by~$T_1 \cup \{m_3\}$.
\item Writing $T_1(k)$ for the clone generated by $T_1 \cup \{ m_k\}$
we have
$$\langle T_1\rangle  \subseteq \cdots
\subsetneqq   T_1(5) \subsetneqq
 T_1(4) \subsetneqq T_1(3) = \Pol(T_1)
\subsetneqq T_1(2) = \O .$$
\item
Every clone in the interval $(\langle
T_1\rangle  ,\O]$ is equal to one of the clones~$T_1(k)$.
\end{enumerate}
\end{Theorem}

Note that this implies that the rank of each clone $T_1(k)$ is~$k-2$,
and the rank of~$\langle T_1\rangle$ is~$\omega$.
\subsubsection{Around the binary clone $\cl{T_2}$}

We now consider the base set $X=\bbn$ only.  Recall that $T_2$ is the
set of all functions $f\in \oo 2 $ which are never 1-1, i.e., not 1-1
even when restricted to any set of the form $(A\times B)\cap \nabla$
or~$(A \times B)\cap \Delta$, with infinite $A$ and~$B$.
This definition seems to be complicated, certainly more complicated
than the definition of~$T_1$.    But there is a reason for this: The clones
$T_2$ and $\Pol(T_2)$ {\em are} complicated (in the sense of
Descriptive Set Theory), as we will see in this section, in particular
in  Theorem~\ref{pi11}.

We will need the following facts and definitions from Descriptive Set
Theory (see \cite{Moschovakis:1980}, \cite{Kechris:1995}):
\begin{itemize}
\item A Polish space is a separable topological space whose topology
      is generated by a complete metric.  Examples
      are~$\bbr$,~$\bbn$, the Cantor space~$2^{\bbn}$, the Baire
      space~$\bbnn$; finite or countable products of Polish spaces are
      again Polish spaces.
\item A subset $Y$ of a Polish space is called \emph{analytic}
 iff $Y$ can be written as the continuous image of a closed subset of~$\bbnn$.  All Borel sets are analytic.
\item The \emph{coanalytic}
 subsets of a Polish space are exactly the complements of analytic sets.
\item If $f: M_1\to M_2$ is a continuous map between Polish spaces, then $f^{-1}[Y] \subseteq M_1$ is analytic (coanalytic) whenever $Y \subseteq M_2$ is analytic (coanalytic).
\item In every uncountable Polish space there are analytic sets which are not coanalytic.
\item A coanalytic set $Y\subseteq M_2$ is called \emph{completely
coanalytic} iff: For every Polish space~$M_1$, every coanalytic subset
$Y' \subseteq M_1$ is the preimage of~$Y$ under some continuous
function~$f:M_1\to M_2$. A complete coanalytic set can therefore not
be analytic.
\end{itemize}

A central theme of Descriptive Set Theory is the  investigation of
the complexity of subsets of Polish spaces.  Borel sets are
considered relatively simple; the simplest of all are of course the
closed and the open sets.  Most sets of real numbers that appear in
analysis are in fact Borel sets.  Analytic sets are more complicated
than Borel sets (similar to the difference between recursively
enumerable
and recursive sets), and coanalytic sets are considered to be slightly more
complicated.

Using a bijection from $\bbn^n$ onto~$\bbn$, each set $\O\un_\bbn =
\bbn ^{\bbn^n} $ can be naturally bijected onto $\bbn^\bbn$ and
becomes thus a Polish space. Also $\O_\bbn = \bigcup_{n\ge 1}
\O_\bbn\un $ can be naturally bijected onto~$\bbn \times \bbn^\bbn$,
which is itself homeomorphic to~$\bbn^\bbn$. We can thus
apply notions from Descriptive Set Theory to sets of operations on
$\bbn$, and in particular measure the complexity of clones.

\begin{Fact} If $\B \subseteq \O$ is a Borel or analytic set, then the clone $\langle \B\rangle$
is analytic.
\end{Fact}

However, in many cases Borel sets of functions will again generate a Borel clone. This motivates the following question:
\begin{prob}
Find a Borel subset $\B \subseteq \O$ (preferably containing $\oo 1 $)
 such that the clone
$\langle \B \rangle$ is not a Borel set.
\end{prob}

\begin{theorem}[{Goldstern} \cite{analytic}]\label{pi11}
The  clone  $\Pol(T_2)$ (as well as the binary Menger algebra $T_2$) is a
complete coanalytic set.\end{theorem}
Contrast this with the following fact about $T_1$:
\begin{theorem}[{Pinsker} \cite{Pin04almostUnary}]
Each of the sets~$T_1(k)$, as well as the sets $T_1$ and~$\langle T_1\rangle$,
are Borel sets.
\end{theorem}

 From Theorem~\ref{pi11} we get:
\begin{cor}\label{pi11.cor} \
\begin{enumerate}
\item Neither $\Pol(T_2)$ nor $T_2$ can be finitely or countably
generated over~$\oo 1 $.
\item The intervals $[\langle \oo 1 \rangle , \Pol(T_2)]$
and even $[\langle \oo 1 \rangle , \langle T_2\rangle]$
are   uncountable.
\end{enumerate}
\end{cor}
\begin{proof}
(1): The set $\oo 1 \cup \C$ is a Borel set for any countable set~$\C$,
   so $\langle \oo 1 \cup \C\rangle$ is analytic, and cannot be equal
   to~$\Pol(T_2)$.

(2): We can find an uncountable sequence $(f_\alpha: \alpha<\omega_1)$
of functions in~$T_2$ such that~$f_\alpha\notin  \C_\alpha$,
where $\C_\alpha $ is the (analytic!) clone generated
by~$ \oo 1 \cup \{f_\gamma: \gamma  < \alpha\}$.  All the clones
$\C_\alpha$
are different.
\end{proof}

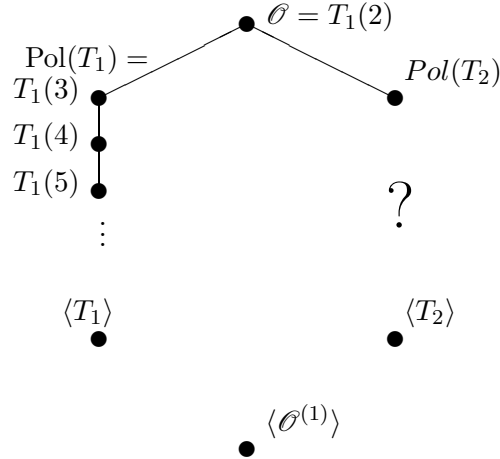
\newcommand{\cto}{\langle T_1\rangle}
\newcommand{\pto}{\Pol(T_1)}
\newcommand\OO{\O}
\def\mm_#1{\hskip-0.4cm T_1(#1)}

\newlength{\normalunitlength}
\setlength{\normalunitlength}{\unitlength}

\setlength{\unitlength}{0.7\unitlength}
    \begin{center}
    \begin{picture}(200,350)

    \put(100,70){\circle*{8}}
    \put(110,80){$\cl{\Oo}$}

    \put(100,300){\circle*{8}}
    \put(110,300){$\OO=T_1(2)$}

    \put(180,260){\circle*{8}}
    \put(185,270){$Pol(T_2)$}

    \put(20,130){\circle*{8}}
    \put(0,140){$\cto$}

    \put(180,130){\circle*{8}}
    \put(185,140){$\langle T_2\rangle$}

    \put(20,260){\line(2,1){80}}
    \put(180,260){\line(-2,1){80}}

    \put(175,190){\Huge{?}}

    \put(20,260){\circle*{8}}
    \put(-20,277){$\pto=$}
    \put(-10,260){$\mm_3$}

    \put(20,210){\line(0,1){50}}

    \put(20,235){\circle*{8}}
    \put(-10,235){$\mm_4$}

    \put(20,210){\circle*{8}}
    \put(-10,210){$\mm_5$}

    \put(20,180){$\vdots$}

    \end{picture}
    \end{center}
\setlength{\unitlength}{\normalunitlength}

We have seen in Theorem~\ref{pinsker:chain} that the interval
$[\langle T_1\rangle , \O]$
is completely understood. The situation is very different with
the clone~$T_2$.

\begin{prob}
\begin{itemize}
\item Is~$\langle T_2 \rangle = \Pol(T_2)$?
\item If not, how many elements does the interval $[\langle T_2\rangle, \Pol(T_2)]$
have?  More generally, what is the lattice-theoretic structure of this
interval?
\end{itemize}
\end{prob}
\subsubsection{Around $\oo 1 $}

We have seen in Theorem~\ref{pi11} that there are uncountably many
clones between $\oo 1 $ and~$\langle T_2\rangle$. The next theorem
shows that the situation is even worse:

\begin{theorem}[{Goldstern, S\'agi and Shelah} \cite{GGS}]
On a countable base set $X$ there are $2^{2^{\aleph_0}}$ clones
containing~$\oo 1 $.  In fact, there is an order-preserving
embedding of the power set of~$\bbr$ into~$[\loor , \cl{T_2}]$.
\end{theorem}

However, the many clones in this theorem appear very low in the
interval~$[\loor  , \langle T_2\rangle]$.
 Note that for every analytic clone $\C\supseteq
\oo 1 $  the interval $[\C , \cl{T_2} ]$ is uncountable (by an
argument similar to the one in Corollary~\ref{pi11.cor}); this
motivates the following question.

\begin{prob}
Assume that $\C \supseteq \oo 1 $ is an analytic clone distinct from
$\O$.  Does $[\C, \cl{T_2}]$ have to have cardinality at
least~$2^{\aleph_0}$?  At least~$2^{2^{\aleph_0}}$?
\end{prob}

\SUBSECTION{The interval of clones that contain all permutations}
    For the case of a base set $X$ of regular cardinality, we have
    seen a complete list of those dual atoms of~$[\cl{\S},\O]$
    which do not contain $\Oo$ (Theorem \ref{thm:heindorfpinsker:maximal:aboveS}).
    Also, we know that the interval $[\cl{\S},\O]$
    is large since its subintervals $[\cl{\Oo},\O]$ (confer the preceding section)
    and $[\cl{\S},\cl{\Oo}]$ (confer the following section) are.
    So an interesting next question is to determine the atoms of
    the interval; analogously to the atoms of~$\Cl(X)$, such
    clones are generated by a single function, which here is not a
    permutation. Now it is a fact that $\pol(\S)=\cl{\S}$ (see
    e.g.\ \cite{MP06MinimalAboveS} for a proof), which has as a
    consequence that all functions $f\in\O$ which generate atoms in
    $[\cl{\S},\O]$ are essentially unary: Indeed, otherwise $f\nin\pol(\S)$, and there would exist
    $g_1,\ldots,g_{n_f}\in\S$ such that the unary operation
    $h=f(g_1,\ldots,g_{n_f})\nin\S$. However, then $\cl{\{h\}\cup\S}$ is a proper subclone
    of~$\cl{\{f\}\cup\S}$ as it contains only essentially unary
    operations, in contradiction with the assumption that $f$
    together with $\S$ generates an atom of~$[\cl{\S},\O]$.
    \begin{ex}
        Any constant operation on~$X$ together with $\S$ generates
        an atom of~$[\cl{\S},\O]$.
    \end{ex}

    \textsc{Pinsker} and \textsc{Machida} \cite{MP06MinimalAboveS} gave an explicit description of all
    operations which together with $\S$ generate atoms of~$[\cl{\S},\O]$, for all
    infinite~$X$. Moreover, they described what these atoms look like; since their theorem is quite technical,
    we do not state it here. As a corollary, they
    found that the number of atoms of~$[\cl{\S},\O]$  on an
    infinite set of cardinality $\aleph_\alpha$ is
    $\max\{|\alpha|,\aleph_0\}$. Applying their result, one finds
    that the atoms of the interval for the countably infinite case
    are the following:

    \begin{itemize}
        \item For every~$n\geq 0$, the clone of all essentially unary operations
            whose corresponding unary operation is either
            a permutation or has only infinite kernel classes and exactly $n$ elements outside its range.
        \item The clone of all essentially unary operations
            whose corresponding unary operation is either a permutation
            or a constant.
        \item The clone of all essentially unary operations whose corresponding unary operation
                is injective and either is a permutation or has an infinite complement of its range.
    \end{itemize}

\SUBSECTION{The interval of unary clones that contain all
permutations}
    Recall that a clone is called unary iff all its operations are
    essentially unary. Unary clones are nothing but disguised
    submonoids of the full transformation monoid~$\Oo$, since they
    arise from such submonoids by adding fictitious variables to the operations of the monoid. The
    interval $[\cl{\S},\cl{\Oo}]$ of unary clones that contain $\S$ is therefore just
    the interval $[\S,\Oo]$ of the lattice of submonoids of~$\Oo$.
    Again by \textsc{Zorn}'s lemma and the fact that $\Oo$ is finitely generated over~$\S$,
    this interval is dually atomic. Its dual
    atoms, which we call \emph{precomplete monoids}, have been described by \textsc{Gavrilov} for countable $X$ and
    by \textsc{Pinsker} for all uncountable~$X$.
    \begin{thm}[Gavrilov \cite{Gav65}, Pinsker
    \cite{Pin05MaxAbovePerm}]\label{thm:gavrilovpinsker:interval:aboveS:maximal}
        If $X$ has regular cardinality, then the
        precomplete submonoids of~$\Oo$ which contain the permutations
        are exactly the monoid $\A$ and the monoids $\G_\lambda$
        and $\M_\lambda$ for $\lambda=1$ and $\aleph_0\leq \lambda\leq |X|$, $\lambda$ a cardinal, where
        \begin{itemize}
            \item{$\A=\{f\in\Oo:|f^{-1}[\{y\}]|<|X|$ for all but fewer than $|X|$ many $y\in X \}$}
            \item{$\G_\lambda=\{f\in\Oo:$ $|X\setminus f[X\setminus A]|\geq \lambda$ for all $A\in [X]^\lambda\}$}
            \item $\M_\lambda=\{f\in\Oo: |X\sm f[X]|<\lambda$ or $f\on {(X\sm A)}$
            is not injective for any $A\in [X]^{< \lambda}\}$.
        \end{itemize}

        If $X$ has singular cardinality, then the same is true
        with the monoid $\A$ replaced by
        \begin{itemize}
            \item
            $\A'=\{f\in\Oo: \exists \lambda < |X| \,\,(\,|f\inv[\{y\}]|\leq\lambda$
            for all but fewer than $|X|$ many $y\in X\,)\,\}$.
        \end{itemize}
    \end{thm}
    \begin{cor}\label{COR:numberOfMaximalMonoids}
        On a set $X$ of infinite cardinality $\aleph_\alpha$ there
        exist $2\,|\alpha|+5$ precomplete submonoids of~$\Oo$ that contain
        the permutations. Hence the smallest cardinality on which there
        are infinitely many such monoids is $\aleph_{\omega}$.
    \end{cor}

    It is interesting to compare the list of precomplete clones above
    $\S$ in Theorem \ref{thm:heindorfpinsker:maximal:aboveS} with
    the one of precomplete monoids above $\S$ in Theorem
    \ref{thm:gavrilovpinsker:interval:aboveS:maximal}: One
    finds that there exist precomplete submonoids of~$\Oo$ whose
    clone of polymorphisms is precomplete, but also precomplete monoids
    whose polymorphism clone is not precomplete, and non-precomplete monoids
    whose polymorphism clone is precomplete.

    Although there exist so few dual atoms in this interval,
    it is
    huge:

    \begin{thm}[Pinsker \cite{Pin05NumberOfUnary}]\label{thm:pinsker:numberOfMonoidsAboveS}
        Let $X$ have cardinality
        $\kappa=\aleph_\alpha$. Then there exist $2^{2^\lambda}$ submonoids of~$\Oo$ which contain
        all permutations, where
        $\lambda=\max\{\,|\alpha|,\aleph_0\}$.\\ Moreover, if $\kappa$
        is regular, then $|[\S,\G]|=2^{2^\lambda}$ for every
        precomplete monoid $\G$ above~$\S$; in fact,
        $|[\S,\D]|=2^{2^\lambda}$, where $\D$ is the 
        intersection of the precomplete elements of~$[\S,\Oo]$.\\ If
        $\kappa$ is singular, then $|[\S,\G]|=2^{2^\lambda}$ for all precomplete
        monoids $\G$
        except $\A'$: If $\lambda<\kappa$, then $|[\S,\A']|=|[\S,\D]|=2^{2^\lambda}$,
        but if $\lambda=\kappa$, then $|[\S,\A']|=|[\S,\D]|=2^{(\kappa^{<\kappa})}$ (where $\kappa^{<\kappa}=\sup\{\kappa^\xi:\xi<\kappa\}$).
    \end{thm}

    In the same paper, it was remarked that if GCH holds, then $2^{(\kappa^{<\kappa})}=2^{2^\kappa}$, so
    in this case we have $|[\S,\D]|=2^{2^\lambda}$ on all
    infinite~$X$. However, for any singular $\kappa$ it is also consistent that $2^\kappa<
    2^{(\kappa^{<\kappa})}<2^{2^{\kappa}}$. Therefore, if $\kappa$ is singular and $\aleph_\kappa=\kappa$,
    then the intervals $[\S,\A']$ and $[\S,\D]$ can be smaller than
    $2^{2^{\lambda}}$. In particular we have that whether or not the intervals $[\S,\A']$
    and, say, $[\S,\M_1]$ are of equal cardinality depends on
    the set-theoretical universe.

    Not only the dual atoms, but also the atoms of the interval
    $[\S,\Oo]$ are known: As we have seen in the preceding
    section, they are just the atoms of the interval
    $[\cl{\S},\O]$ of the clone lattice listed in the article \cite{MP06MinimalAboveS}, since all such atoms
    turned out to be unary.

\SUBSECTION{The interval of clones above the idempotent
clone}\label{sec:aboveIdempotent}

Consider the clone $\{f:  \forall x\, \, f(x,\ldots, x) = x\}$,
which consists of all idempotent functions. For any subset $A
\subseteq X$ the set
$$ \{ f: \forall x\in A\,\, f(x,\ldots, x) = x\}$$
is again a (larger) clone. More generally, we use the following
definition:

\begin{Definition}
Let ${D}$ be a filter on~$X$, that is: ${D}$
is a nonempty subset of the power set of~$X$ which is upward closed
and closed under finite
intersections. We allow here also the improper filter consisting of
all subsets of~$X$.

Then we define
$$ \C_{{D}}:=
\{ f: (\exists A\in {D})(\forall x\in A) \, f(x,\ldots, x) = x\}.$$
\end{Definition}

Observe that with this definition the clone of idempotent operations
is just $\C_{\{X\}}$. The following theorem gives an example of a
rather complicated interval in the clone lattice, whose structure is
nevertheless ``known''.

\begin{theorem}[{Goldstern and Shelah} \cite{747}]
\ \label{filterclones}
\begin{itemize}
\item
    For every filter ${D}$ the set $\C_{D}$ is a clone.
\item For ${D}= \{X\}$, $\C_{D}$ is the clone of all idempotent functions,
    and if ${D}$ is the improper filter (containing $\emptyset$), then
    $\C_{D}=\O$.
\item
   For filters ${D}_1 \subsetneqq {D}_2$ we have $\C_{{D}_1} \subsetneqq \C_{{D}_2}$.
\item
   Every clone in the interval $[\C_{\{X\}},\O]$ is equal to some~$\C_{D}$.
\end{itemize}
\end{theorem}


In particular, each ultrafilter on~$X$ corresponds to a precomplete
clone on~$X$.  Note that unlike the map $I\mapsto \C_I$ defined in
Section~\ref{sec:maximal:number}, the map $D \to \C_D$ is monotone.

We write $\beta X$ for the set of ultrafilters on~$X$.
Recall that $\beta X$
 carries a topological structure: For each $A \subseteq X$ the
set $\hat A:= \{U\in \beta X:  A\in U\}$ is declared open, and these
sets $\hat A$ form a basis.   $\beta X$ is also known as
 the \textsc{Stone}-\textsc{\v{C}ech}
 compactification of the discrete space~$X$, see \cite{Comfort-Negrepontis}.

The filters ${D}$ on~$X$ are in natural 1-1 correspondence with the
closed subsets of~$\beta X$ through the map$$ {D} \mapsto \{ U: {D}
\subseteq U\}. $$ We thus obtain an isomorphism between the lattice
of closed subsets of~$\beta X$ (ordered by reverse inclusion) and
the lattice $[\C_{\{X\}}, \O]$: The empty set corresponds to the
full clone~$\O$, the points in~$\beta X$ correspond to the
precomplete clones, and larger sets correspond to smaller clones.
This correspondence, even though it is quite straightforward, is an
example to show that questions about the structure of the clone
lattice can sometimes be translated to questions about better known
topological spaces.  For example, using the following fact from
topology, we can  compute the rank of any clone in the interval
$[\C_{\{X\}},\O]$.   Ranks of clones are translated to ranks of
closed subsets of~$\beta X$; the rank of a closed set $C$ is equal
to $\alpha$ if it is not equal to any $\beta<\alpha$, but every
closed set strictly contained in~$C$ has rank~$<\alpha$.

The following fact is easy to prove:
\begin{Fact}
Let $X$ be a discrete space, and let $A \subseteq \beta X$ be a closed infinite set.  Then $A$ contains a homeomorphic copy of~$\beta\bbn$.
\end{Fact}
\begin{cor} Let $\C\in [\C_{\{X\}},\O]$ be a clone.
\begin{enumerate}
\item
If $\C$ is the intersection $\C_{D_1}\cap \cdots \cap \C_{D_n}$, where
the $D_i$ are distinct ultrafilters, then $\rank(\C)=n$.
\item Otherwise, $\rank(\C)=\infty$.
\end{enumerate}
\end{cor}

\begin{proof}
(1): The rank of any finite subset of~$\beta X$ is equal to its
cardinality.

(2): The space~$\beta\bbn $ contains a strictly decreasing sequence of
closed sets: $\beta \bbn \supsetneqq \beta(\bbn \setminus
\{0\})\supsetneqq \cdots $.
\end{proof}

\subsubsection{Intervals below sufficiently rich clones}
Since we know that there are always $2^{2^{|X|}}$ ultrafilters on any
 infinite set~$X$, Theorem~\ref{filterclones} gives an easy proof of
 the \textsc{Gavrilov}-\textsc{Rosenberg} theorem from \cite{Gav65}
 and \cite{Ros76} that $\Cl(X)$ has the largest possible number of
 coatoms.  This particular result can also be obtained from the
 special case $\C=\O$ of
\textsc{Marchenkov}'s Theorem~\ref{marchenkov} below.

\begin{Definition}  For any clone $\C$ let $B(\C):= \{\nu(f): f\in \C\}$,
where for any function $f\in \O $ we define $\nu(f)= \{x\in X:
f(x,\ldots, x)=x\}$.
\end{Definition}
\begin{definition}
The \emph{discriminator function} $d\in \oo 3 $ is defined by
$$ d(x,y,z) =
 \begin{cases}
x & \mbox{if $x=y$}\\
z & \mbox{otherwise.}
\end{cases}
$$
\end{definition}

\begin{Theorem}[Marchenkov
 \cite{Mar81-CardinalityOfTheSetOfPrecomplete}]
\label{marchenkov} Let $\C$ be a clone such that $B(\C)$ is a Boolean
subalgebra of the power set of~$X$, and assume~$d\in \C$.  Then every
ultrafilter $U$ on the Boolean algebra $B(\C)$ induces a clone$$\C_U:=\{ f\in \C: \nu(f)\in U \}$$
which is a coatom in the interval~$[\J, \C]$.
\end{Theorem}

\SUBSECTION{Monoidal
intervals}\label{sec:intervals:monoidalIntervals}
    Let $\M\subseteq\Oo$ be a submonoid of the full transformation
    monoid~$\Oo$. Then the set of those clones $\C$ which have $\M$
    as their 
    unary fragment (i.e., which satisfy $\C\uo=\M$) is an
    interval of the clone lattice: Clearly, the smallest clone with
    this property is the clone $\cl{\M}$ which consists of all
    essentially unary functions whose corresponding unary function
    is an element of~$\M$. The top of the interval is easily seen to
    be the set $\pol(\M)$ of all functions that preserve the
    monoid~$\M$. 
    \begin{ex}
        The monoidal interval corresponding to the full
        transformation monoid $\Oo$ is just the interval
        $[\cl{\Oo},\O]$, and has been subject to much investigation; see Section \ref{sec:interval:aboveOo}.
        In particular, we know that if $X$ is countably infinite,
        then the cardinality of this interval
        equals $|\Cl(X)|=2^{2^{\aleph_0}}$.
    \end{ex}
    \begin{ex}
        The monoid $\S$ of permutations of~$X$ has a monoidal interval which consists of just one
        element: $\cl{\S}=\pol(\S)$ (see e.g.\
        \cite{MP06MinimalAboveS}). Monoids, and also clones, whose monoidal interval
        has this property are called \emph{collapsing}.
    \end{ex}
    \begin{ex}
        Those clones which contain only idempotent operations form
        exactly the monoidal interval induced by the trivial monoid
        $\{\pi^1_1\}$, since $\pol(\{\pi^1_1\})$ is the clone of
        idempotent functions.
    \end{ex}

    One reason why monoidal intervals are interesting is that
    studying such an interval is in some sense ``orthogonal'' to
    studying the lattice of
    monoids: In the first case, we fix the monoid, and look how
    functions of larger arity generate each other modulo that
    monoid, whereas in the latter case we forget about higher
    arities and concentrate on unary operations only. This way,
    because the monoidal intervals are a partition of~$\Cl(X)$,
    investigating the clone lattice is split into the study of
    monoidal intervals and the study of the monoid lattice.

\begin{figure}
    \begin{picture}(210,210)(-105,0)
        {\qbezier(0,5)(-90,40)(-100,150)}
        {\qbezier(0,5)(90,40) (100,150)}
        {\qbezier(0,5)(-50,50)(0,110)}
        {\qbezier(0,5)(50,50) (0,110)}
        \qbezier(0,5)(-50,90)(-86,90)
        \qbezier(0,5)(-60,35)(-86,90)
        \put(0,5){\circle*{8}}
        \put(-3,16){$\J$}

        \put(-86,90){\circle*{8}}
        \put(-86,97){$\pol(\{\pi^1_1\})$}

        \qbezier(-12,70)(-50,100)(-60,189)
        \qbezier(-12,70)(-10,100)(-60,189)
        \put(-12,70){\circle*{8}}
        \put(-12,56){$\cl{\M}$}

        \put(-60,189){\circle*{8}}
        \put(-69,197){$\pol(\M)$}

        \put(0,110){\circle*{8}}
        \put(6,110){$\cl{\Oo}$}
        \qbezier(0,110)(-35,150)(0,200)
        \qbezier(0,110)(35,150) (0,200)

        \put(0,200){\circle*{8}}
        \put(10,203){$\O=\pol(\Oo)$}
        {\qbezier(0,200)(-90,195)(-100,150)}
        {\qbezier(0,200)(90,195) (100,150)}
    \end{picture}
\end{figure}

    There is another concept justifying the study of monoidal
    intervals. For two distinct clones $\C$ and~$\D$, there exists $n\geq
    1$ such that $\C\un\neq\D\un$. Moreover,
    if this is the case and~$k\geq n$, then also $\C\uk\neq\D\uk$.
    Therefore, one could say that $\C$ and $\D$ are closer the later
    their $n$-ary fragments start to differ. More precisely,
    $$
        d(\C,\D)=\begin{cases}\frac{1}{2^{n-1}}&\C\neq\D\wedge
        n=\min\{k:\C\uk\neq\D\uk\}\\0&\C=\D\end{cases}
    $$
    defines a metric on the clone lattice, first introduced by
    \textsc{Machida} \cite{Mac98} (for a finite base set, but the same works on infinite sets).
    Formulated in this metric, a
    monoidal interval is just an open ball of radius $1$ in the metric
    space $(\Cl(X),d)$. It also makes sense to consider refinements of this
    partition, for example open balls of radius $\frac{1}{2}$, or
    equivalently sets of clones with identical binary fragments;
    they are of the form $[\cl{\h},\pol(\h)]$, where $\h\subseteq\Ot$
    is a binary Menger algebra.
    

    On finite~$X$, monoidal intervals are either finite, countably
    infinite, or of size continuum (first mentioned in \cite{RSxxIntervalCardinality}, see also
    the introduction of \cite{Pin06Monoidal} for a proof using Descriptive Set Theory). Moreover, those
    possibilities are all realized: There must be monoidal
    intervals of size continuum for cardinality reasons, the monoidal
    interval corresponding to~$\Oo$ is a finite chain of length
    $|X|+1$ \cite{Bur67}, and a countably infinite monoidal interval was exposed
    by \textsc{Krokhin} \cite{Kro97}. The question whether on infinite $X$ monoidal intervals can have ``strange''
    cardinalities, i.e., cardinalities
    strictly between $|X|$ and~$2^{|X|}$, has recently been
    answered:

    \begin{thm}[Pinsker
    \cite{Pin06Monoidal}]\label{thm:pinsker:monoidal:cardinality}
        There exist at least monoidal intervals
        of the following cardinalities:
        \begin{itemize}
            \item{$\lambda$ for all $\lambda\leq 2^{|X|}$.}
            \item{$2^\lambda$ for all $\lambda\leq 2^{|X|}$.}
        \end{itemize}
    \end{thm}

    One might ask if not all cardinals $\le 2^{2^{|X|}}$ can appear
    as cardinalities of monoidal intervals; but this is consistently false,
    by the following consequence of a theorem due to \textsc{Kunen} which has been pointed out by
    \textsc{Abraham}. He observed that the following conditions
    (a), (b), and (c) are consistent:

        \begin{itemize}
        \itm a  $2^{\aleph_0} = \aleph_1$
        \itm b $2^{\aleph_1} > \aleph_2$
        (in fact,  $2^{\aleph_1} $ can be arbitrarily large)
        \itm c Whenever $\F$ is a family of subsets of~$\omega_1$ which is
        closed unter arbitrary intersections and arbitrary increasing
        unions, then $\F$ has either $\le \aleph_1$ elements, or $\ge
        2^{\aleph_1}$ elements.
        \end{itemize}

    Now (a), (b) and (c) together imply
    \begin{itemize}
        \itm d All intervals in the clone lattice on
        a countably infinite base set have cardinality $\le \aleph_1$ or
        $= 2^{\aleph_1}$. In particular, there is no interval of
        cardinality~$\aleph_2$.
    \end{itemize}

    However, this is really a remark about cardinalities of algebraic
    lattices, since (a), (b) and (c) together also imply

    \begin{itemize}
        \itm e Every algebraic lattice with $\aleph_1$ compact
        elements has cardinality $\le \aleph_1$ or
        $= 2^{\aleph_1}$. In particular, there is no algebraic
        lattice with $\aleph_1$ compact elements that has
        cardinality~$\aleph_2$.
    \end{itemize}

    Therefore, the ``right'' question to ask is the following:

    \begin{prob}
        Is every algebraic lattice with at most $2^{|X|}$ compact elements equipotent to a
        monoidal interval of the clone lattice?
    \end{prob}

    In case of a negative answer, the same can be asked about arbitrary, not necessarily
    monoidal intervals, which leads to the following
    easier variant of Problem
    \ref{prob:algebraiclattice:IsomorphicToInterval}:

    \begin{prob}
        Is every algebraic lattice with at most $2^{|X|}$ compact elements equipotent to a
        an interval of the clone lattice?
    \end{prob}

    Theorem \ref{thm:pinsker:monoidal:cardinality} was in fact a
    corollary of a result in the same paper on the possible
    structure of monoidal intervals:

    \begin{thm}[Pinsker
    \cite{Pin06Monoidal}]\label{thm:pinsker:monoidal:structure}
        Let $\L$ be an algebraic and dually algebraic distributive lattice with at
        most $2^{|X|}$ completely join irreducible elements. Then there is a monoidal interval in
        $\Cl(X)$ isomorphic to~$0+\L$, which denotes $\L$ with an additional smallest element $0$ added.
    \end{thm}

    We remark here that the class of algebraic and dually algebraic distributive
    lattices is the class of completely distributive
    lattices, or equivalently the class of lattices of order ideals of
    partial orders (see e.g.\ \cite[p.83]{CD73} for the latter statement).

    It is not surprising that the monoidal intervals of this theorem
    are not all possibilities: For example, it has been remarked in the same paper that
    the monoidal interval of idempotent clones is not modular.

    \begin{prob}
        Find other classes of lattices that appear as monoidal
        intervals.
    \end{prob}

    It could even be the case that all algebraic lattices which
    satisfy the only obvious restriction of not having more than
    $2^{|X|}$ compact elements appear as monoidal intervals.

    \begin{prob}
        Is there an algebraic lattice with at most
        $2^{|X|}$ compact elements that is not isomorphic to a monoidal interval?
    \end{prob}

\SECTION{The local clone lattice}\label{sec:local}

    Fix some index set~$I$, and let $R\subseteq X^I$ be an $I$-ary relation,
    i.e., a set of~$I$-tuples with entries in~$X$. If $f\in\O$,
    then we say that $f$ \emph{preserves} $R$ iff
    $f(r_1,\ldots,r_{n_f})\in R$ for all $r_1,\ldots,r_{n_f}\in
    R$ ($f(r_1,\ldots,r_{n_f})$ denotes the $I$-tuple that results if we apply $f$ to the tuples $r_j$
    componentwise; this notation also agrees with the ``composition'' notation
    introduced at the beginning  of Section~\ref{sec:1}). We have seen this concept earlier in this
    paper: If $I={X^n}$ for some natural number $n\geq 1$, then
    $R\subseteq X^{X^n}$ is just a set of~$n$-ary operations and $f$ preserves $R$
    iff $f\in\pol(R)$. Another important case is when $I$ is
    a positive natural number, and $R$ is a \emph{finitary relation} on~$X$.

    Now we set $\pol(R)\subseteq\O$ to consist of all operations that
    preserve~$R$, for an arbitrary relation~$R$; this definition is an extension of the case
    where $R$ is a set of operations. We write
    $\ppol(\R)=\bigcap\{\pol(R):R\in\R\}$ for a \emph{set of relations}~$\R$. Conversely, for a set of operations
    $\F\subseteq\O$, we write $\iinv(\F)$ for the set of all
    finitary relations that are preserved by all~$f\in\F$.

    On a finite base set, the operators $\ppol$ and $\iinv$ are a
    strong tool for describing clones, since in that case finitary
    relations suffice to describe all clones: If $\C\subseteq\O$
    is any set of operations, then $\C$ is a clone iff
    $\C=\ppol\iinv(\C)$. In other words, every clone $\C$ is
    determined by the finitary relations it preserves. Moreover,
    for any $\C\subseteq\O$, $\cl{\C}=\ppol\iinv(\C)$. Since clones preserve fewer relations
    the larger
    they are, this method is particularly useful when describing large
    clones, such as precomplete ones, whereas small clones are often
    better described by (functional) generating systems.

    If $X$ is infinite, then not every clone is of the form
    $\ppol(\R)$ for a set $\R$ of finitary relations, although sets
    of operations of this form are still clones.
    \textsc{Rosenberg} \cite{Ros71-AClassOfUnivAlgByInfinitaryRel} observed first that every clone is
    of the form $\ppol(\R)$ for a set of infinitary relations: For example, it follows from our observations
    in the beginning of Section \ref{sec:maximal} that
    $\{\C\uo,\C\ut,\ldots\}$ is such a set for any clone~$\C$.

    A clone is called \emph{locally closed} or \emph{local} iff it
    is of the form $\ppol(\R)$ for some set of finitary
    relations~$\R$. This naming is made clearer by the following
    alternative 
    definition which is easily seen to be equivalent. A clone is
    local iff it satisfies the following additional closure
    property: For every~$f\in\O$, if for all finite $A\subseteq X$
    there exists $g\in \C$ of the same arity as $f$ which agrees with $f$
    on~$A$, then~$f\in\C$. This is the same as to say that
    for every $n\geq 1$, $\C\un$ is closed in the product topology (Tychonoff topology) on~$X^{X^n}$,
    where $X$ is
    taken to be discrete. In other words, $\C$ is local iff it
    contains all operations that can be ``locally''
    approximated (i.e., interpolated on every finite set)
    by functions from~$\C$. This
    generalization of the Galois connection $\iinv-\ppol$ from the case
    where $X$ is finite is due to \textsc{Romov} \cite{Rom77}.

    To emphasize the distinction between clones and local clones,
    we may call the former \emph{global} clones; that is, a global
    clone is just a clone which is not necessarily local.
    Arbitrary intersections of local clones yield local clones, and
    the local clones on~$X$ form a complete lattice $\Cl_{loc}(X)$,
    which is \emph{not} a sublattice of~$\Cl(X)$: In general, the
    clone generated by two local clones (in $\Cl(X)$) need not be locally closed.

    \begin{ex}
        Let $X$ be the set of integers $\mathbb{Z}$, and let
        $f \in\S$ be the permutation that switches $0$ and~$1$,
        and is the identity otherwise. Let $g\in\S$ be the
        permutation that maps every $x\in \mathbb{Z}$ to~$x+1$.
        Consider the local clones $\C$ and $\D$ locally generated by
        $\{f\}$ and $\{g,g\inv\}$, respectively. The only non-trivial
        unary operation in~$\C$ is~$f$, and in~$\D$ we only
        have the operations $g^k$ and~$g^{-k}$, for all
        $k\geq 1$. Now one can verify that the join of~$\C$ and $\D$
        in~$\Cl_{loc}(X)$ contains~$\S$, which cannot be the case
        for the join in the global clone lattice~$\Cl(X)$, since this join is
        countable but $|\S|=2^{\aleph_0}$.
    \end{ex}

We are grateful to the referee for pointing out that $\Cl_{loc}(X)$ is not
algebraic.   In fact, $\Cl_{loc}(X)$ has no compact elements, except for 
the clone of all projections: 

\begin{ex}
Let $X$ be infinite, and fix a linear order $\le$ on~$X$
without last element. 
   For each $a\in X$ let 
$$
\begin{array}{rl}
\C_a &:= \{f\in \O: \forall \vec x\,   f(\vec x)  \le a \} \\
\D_a &:= \{f\in \O: \forall \vec x\, ( \max(\vec x)\ge a 
\ \Rightarrow  \ f(\vec x) \ge \max (\vec x) \, \} 
\end{array}
$$
Then 
\begin{enumerate}
\item $\cl{\C_a} = \C_a \cup \{\,\mbox{projections}\,\}$.  
$\cl{\C_a}$ is a local clone.
\item $\cl{\D_a} $ is the set of all functions which are essentially in 
$\D_a$ (i.e., except for dummy variables).
$\cl{\D_a}$ is also a local clone. 
\item If $ a \le a'$, then $\C_a \subseteq \C_{a'} $ and $\D_a \subseteq \D_{a'}$, hence every finite union of clones $\C_a$ (or~$\D_a$, respectively)
is again a clone of this form. 
\item The local closure of~$\bigcup_a\cl{ \D_a}$, as well as the local 
closure of~$\bigcup_a\cl{ \C_a}$, is the clone of all functions. 
\item If $f\in \O$ has unbounded range, 
then $f\notin  \bigcup_a \cl{\C_a}$ (unless $f$ is a projection).
\item If $f\in \O$ has bounded range, then $f\notin  \bigcup_a \cl{\D_a}$. 
\item No local clone  $\C$ (other than the clone of projections)  is 
compact in~$\Cl_{loc}(X)$;  if $\C$ contains a nontrivial unbounded 
function, this is witnessed by the family $(\C_a: a\in X)$, and if 
$\C$ contains a bounded function this is witnessed by the family 
$(\D_a: a\in X)$. 
\end{enumerate}
We leave the proof to the reader. 
\end{ex}

    $\Cl_{loc}(X)$
    is not dually atomic, an example of which was provided
    by \textsc{Rosenberg} and \textsc{Schweigert}
    \cite{RS82-LocallyMaximalClones}, using a relational approach.
    They essentially anticipated the following example:

    \begin{ex}
        For all $n\geq 2$, set $\K_n$ to consist
        of all operations on~$X$ which are either essentially unary or which
        take at most $n$ values. We call $f\in \On$
        \emph{quasilinear} iff there exist functions $\phi_0: 2\To X$ and $\phi_1,\ldots,\phi_n:X\To 2$ such that
        $f(x_1,\ldots,x_n)=\phi_0(\phi_1(x_1)\dot{+}\ldots \dot{+}\phi_n(x_n))$, where $\dot{+}$ denotes the sum
        modulo~$2$. We write $\B$ for the clone of all operations which are either essentially unary or quasilinear;
        $\B$ is often referred to as \textsc{Burle}'s clone.
        Then the interval of non-trivial local clones which contain
        $\Oo$ is
        the following countably infinite chain which ascends to~$\O$:
        $$
        \cl{\Oo}\subsetneqq {\frak B} \subsetneqq \K_2\subsetneqq
\K_3\subsetneqq
        \ldots \subsetneqq \O
        $$
    \end{ex}

    We remark that for finite~$X$, the interval of clones above
    $\Oo$ is exactly this chain, but stops at $\K_{|X|}=\O$
    \cite{Bur67}. Our example then follows easily from the latter fact and the use of local
    closure. As we have seen before, the interval of global clones
    which contain $\Oo$ is fairly complicated, so this
    example supports the intuition that the local clone lattice is
    closer to the clone lattice on a finite base set than to the global clone lattice
    on an infinite base set.

    This idea is also suggested by the fact that the  number of
    local clones is~$2^{|X|}$, so in particular on countably
    infinite $X$ there are as many local clones as there are clones
    on the three-element set. To prove that there are not more local
    clones, it is sufficient to see that a local clone is determined
    by all restrictions of its operations to finite subsets of~$X$;
    for such a set of restrictions, there are not more than
    $2^{|X|}$ possibilities.

    Unfortunately, knowledge of the dual atoms in~$\Cl_{loc}(X)$ is
    not sufficient for a local completeness criterion, i.e., a criterion which says when the local clone generated by
    an algebra equals~$\O$, since $\Cl_{loc}(X)$ is not dually atomic. \textsc{Rosenberg} and
    \textsc{Szab\'{o}} \cite{RS84} gave an example of a cofinal set in
    $\Cl_{loc}(X)$, i.e., a set $\Theta\subseteq\Cl_{loc}(X)$ such that
    every non-trivial local clone is contained in one of the clones
    of~$\Theta$. An algebra is therefore locally complete (that is, it
    locally generates all operations) 
    if and only if its functions are not contained in any of the
    clones of~$\Theta$. Of this system, some clones are dual atoms of~$\Cl_{loc}(X)$ and some
    are not; it does not provide a complete list of the dual atoms. There has been an improvement upon this result in
    \cite{Ros00-LocallyMaximalClonesII}.

    \begin{prob}
        Find all dual atoms of~$\Cl_{loc}(X)$.
    \end{prob}

    The number of dual atoms in~$\Cl_{loc}(X)$ is~$2^{|X|}$,
    since $\pol(A)$ is precomplete (even in~$\Cl(X)$) for every subset $A$ of~$X$, and
    since there do not exist more local clones than that.

    \textsc{Goldstern} and \textsc{Shelah} \cite{747} gave
    examples of fairly complicated intervals in the local clone
    lattice:

    \begin{thm}[Goldstern and Shelah \cite{747}]
        Let ${\frak S}$ be an arbitrary semilattice on~$X$. Then
        there exists a local clone $\C$ such that the interval
        $[\C,\O]$ of the local clone lattice is isomorphic to the
        congruence lattice of~${\frak S}$.
    \end{thm}

    In that paper, the authors remarked that the theorem implies that there exists a local clone $\C$ such that
    $[\C,\O]$ is antiisomorphic to the power set of the natural numbers.
    They also exhibited the following interval which is rather complicated, but
    can in some way be described.

    \begin{ex}
        Let $s$ be a permutation on~$X$ which has only infinite
        cycles, and denote the graph (as a subset of~$X^2$) of~$s$
        by~$s^\circ$. Then the interval $[\pol(s^\circ),\O]$ is 
        isomorphic to the natural numbers $\mathbb{N}$ ordered by the
        divisibility relation, where $1$ is the smallest and $0$ is
        the greatest element.
    \end{ex}

    Observe that this theorem is another example showing that the
    local clone lattice is not dually atomic.

    \begin{prob}
        Which lattices embed into the lattice of local clones?
    \end{prob}

    The idea of local approximation (as interpolation on small sets)
 can be generalized as follows:
    If $X$ is uncountable and $\lambda\leq |X|$, then one can
    define a clone to be $\lambda$-locally closed iff it contains all
    operations which can be interpolated by functions from the clone
    on sets of size smaller than~$\lambda$. With this definition,
    $\aleph_0$-locally closed clones are exactly the local clones.
    It turns out that such clones are precisely the
    polymorphism clones of sets of relations of arity less
 than~$\lambda$, as has been observed in 
    \cite{RS82-LocallyMaximalClones}.
\SECTION{Abstract clones}

\newcommand{\Cf}{{\mathfrak C}}
\newcommand{\Df}{{\mathfrak D}}

\begin{definition}
An \emph{abstract clone} is a many-sorted algebra $\Cf $ consisting
of
\begin{enumerate}
\item infinitely many sorts, i.e., disjoint sets~$C\un$, $n=1,2,\ldots$;
\item operations $*_k^n$ for all $n,k\in \{1,2,\ldots\}$, where
    $*_k^n$ is a map from $C\uk\times (C\un)^k$ to~$C\un$ \\
(we will write $f*(g_1,\ldots, g_k)$ for $*_k^n(f, g_1,\ldots,
g_k)$);
\item constants $p^n_k\in C\un$, for $1\le k\le n$;
\end{enumerate}
and satisfying the following natural set of identities:
\begin{itemize}
\item  $f*(p^n_1, \ldots,p^n_n) = f$ for all $f\in C\un$;
\item $p^n_i*(g_1, \ldots, g_n) = g_i$, whenever $g_1,\ldots, g_n$ are
        in the same sort;
\item the associativity law:
Whenever
\begin{enumerate}
\item[-]  $f\in C\uk$,
\item[-]  $\vec g = (g_1,\ldots, g_k)$, where all $g_i\in C\uell$,
\item[-] $\vec h = (h_1,\ldots, h_\ell)$, where all $h_j\in C\um$,
\end{enumerate}
then
$$ f* \bigl(\,  g_1*(\vec h), \ldots, g_k*(\vec h) \, \bigr) =
 \bigl( f*(g_1,\ldots, g_k)\bigr)*(\vec h).$$
\end{itemize}
\end{definition}

The clones that we have discussed in previous sections will now be
called {\em concrete clones}.  Clearly every concrete clone can be
viewed as an abstract clone: The projections $\pi^n_k$ are
  the constants~$p^n_k$, and $f*(g_1,\ldots, g_k)$ is the functional
composition $f(g_1,\ldots, g_k)$.

\begin{nota}
Whenever $\C$ is a concrete clone, we write $\mathfrak C$ for the corresponding abstract clone (i.e., the multisorted algebra whose universe is the set $\C$),
similarly for the pairs $\C_1/\mathfrak C_1$, $\D/\mathfrak D$, etc.
\end{nota}

\begin{ex}
  Let $\C$ be the (unique) clone on a set with one element.
Then in  the abstract clone $\mathfrak C$, the $n$-th sort is a
singleton with unique element $p^n_1=\cdots = p^n_n$.\end{ex}
\begin{ex}
 Let $\J$ be the set of projections on a set $X$ with
at least two elements. Then (independently of the base set) $\mathfrak J$
is an abstract clone whose only elements are the objects~$p^n_k$, which
are all distinct.
\end{ex}
\begin{ex}
 Let $(L_1, \vee_1)$ and $(L_2,\vee_2)$ be two semilattices, each with at
least two elements. Write $\C_1$ and $\C_2$ for the term clones of
the respective semilattices. Then $\C_1$ and $\C_2$ are {\em abstractly
isomorphic}, i.e., the abstract clones $\mathfrak C_1$ and
$\mathfrak C_2$ are isomorphic as multisorted universal algebras.\end{ex}

The following folklore theorem shows that all abstract clones can be realized
as concrete clones:
\begin{Theorem}[Cayley's theorem]\label{cayley}
Let $\Cf$ be an abstract clone.  Then there is a set $X$ and a concrete clone
$\D $ on~$X$ such that $\Cf $ is isomorphic to~$\mathfrak D$.

In fact, $\D$ can be chosen to be a local clone.
\end{Theorem}
\begin{proof}[Proof sketch]
For $k\le n$ the map $G^n_k: C\uk\to C\un$, defined by~$G^n_k(f)=
f*(p^n_1,\ldots, p^n_k)$, naturally embeds $C\uk$ into~$C\un$.  (The
map corresponds to the operation of adding dummy variables $x_{k+1},
\ldots , x_n$ to make a $k$-ary function into an $n$-ary function.)
Note that $G^n_\ell\circ G^\ell_k = G^n_k$ whenever $k\le \ell\le n$.
Define an equivalence relation $\sim$ as the symmetric closure of the
relation $\{(f,g): \exists n,k\,\, f=G^n_k(g)\}$.
 Let $X:= C/\ssim$.
Note that every equivalence class will have representatives in all except finitely many sorts.

For each $f\in C\un$ we now define a function $\bar f\in\On_X$ as
follows:
 $$
 \bar f( g_1/\ssim, \ldots , g_n/\ssim) = \bigl(f*(g_1,\ldots,
 g_n)\bigr)/\ssim$$
whenever $g_1,\ldots, g_n$ are in the same sort.  It is easy to
check that this definition is proper, that
 $\D := \{\bar f: f\in C\}$ is a local clone (with $\bar p^n_k = \pi^n_k$),
and that the map $f\mapsto \bar f$ is an isomorphism of abstract clones.
\end{proof}

A fundamental  question in the investigation of abstract clones is
the following: Which properties of a concrete clone $\C$ can be read
off of its abstract version $\mathfrak C$?  The examples above show
that the cardinality of the base set is in general not encoded in
the abstract clone, and also the question whether $\C$ is local
cannot be found out by looking at~$\mathfrak C$ only.

However, the question about locality can become interesting when we also
restrict the cardinality of the base set, as the following example shows:

\begin{ex} There is a clone $\C$ on the countable base set $\bbn$ such that
$\mathfrak C$ is not isomorphic to~$\mathfrak D$, whenever $\D$ is a local clone on a countable base set.
\end{ex}

\begin{proof}
For notational simplicity we consider binary Menger algebras instead
of clones. 

We claim that every (binary) concrete clone $\C$ on~$\bbn$ containing
all constant functions which is abstractly isomorphic to a local clone
on $\bbn$ must be an analytic set (and we have already seen that there
are clones on~$\bbn$ which are not analytic).  We use the fact (see
Section~\ref{sec:local}) that local clones are closed (hence
analytic).
Assume that $\D$ is a local
clone on~$\bbn$ corresponding to the abstract clone~$\mathfrak D$,
 and assume that $\iota:\mathfrak C\to \mathfrak D$ is
an isomorphism between abstract clones.  For each $n$ let $c_n$ be the
unary function which is constant with value~$n$, and let $d_n:=
\iota(c_n)$.  It is easy to see that the sets
$$ A := \{(F,G) \in \oo 2 \times \oo 2 : \forall n\, \forall k \,
\forall i \ G(d_n(i),d_k(i)) = d_{F(n,k)}(i)\}$$ as well as~$A \cap
(\oo 2 \times \D)$ are closed, and that $\C$ is just the first
projection of~$A$, hence (as the continuous image of a closed set)
analytic.
\end{proof}

\begin{ex}
Let $\C$ be a precomplete clone on a base set~$X$. $\C$ must contain
$2^{|X|}$ many binary functions.  Consider the abstract clone
$\mathfrak C$; Cayley's construction in Theorem~\ref{cayley} will find
a concrete clone $\D$ on a base set $Y$ of cardinality $|Y| = 2^{|X|}$
such that $\mathfrak C \simeq \mathfrak D$. By a cardinality argument,
$\D$ cannot be precomplete on~$Y$.\end{ex}

Hence, being precomplete is not an ``abstract'' property, but being minimal is:
\begin{Fact}\label{abstract:minimal}
Let $\C$ be a minimal clone.  Then every clone $\D$ with
$\Df \simeq \Cf$ is also minimal.
\end{Fact}
\begin{proof}  Let $\C = \langle\{f\}\rangle$. Let $\Sigma $ be a list of
all equations (in the language of abstract clones, using the
operations $*^n_k$ and the constants $p^n_k$) that $f$ satisfies.  Let
$\iota: \Cf\to \Df$ be an isomorphism; then also $\iota(f)$ will
satisfy the same equations; they will witness that every function
$g\in \langle \{f\}\rangle$ that is not a projection again generates
$f$, and thus satisfies
$\langle \{g\}\rangle =  \langle \{f\}\rangle$.
\end{proof}

\begin{prob}\label{again:prob:minimal:infinitelyManyBinary}
Is there an abstract clone $\Cf$ which is minimal such that $\Cf\ut$
is infinite?
\end{prob}
Note that this problem is really a rephrasing of
 Problem~\ref{prob:minimal:infinitelyManyBinary}.

\markboth{Bibliography}{Bibliography}

\end{document}